\documentclass{article}
\input cyracc.def

\usepackage{amsmath}
\usepackage{amssymb}
\usepackage{amsfonts}
\usepackage{eucal}
\usepackage{url}
\usepackage{bm}
\usepackage[dvipdfmx]{hyperref}
\usepackage{mathrsfs}
\usepackage{graphicx}
\usepackage{amscd}
\usepackage{amsthm}
\usepackage[all]{xy}
\usepackage{tikz}
\usetikzlibrary{cd}
\usetikzlibrary{decorations.pathmorphing}
\usetikzlibrary{patterns}
\DeclareMathOperator{\Ad}{Ad}
\DeclareMathOperator{\Aut}{Aut}
\DeclareMathOperator{\coker}{coker}

\DeclareMathOperator{\ev}{ev}

\DeclareMathOperator{\Hom}{Hom}
\DeclareMathOperator{\id}{id}
\DeclareMathOperator{\Image}{Im}

\DeclareMathOperator{\Ker}{Ker}
\DeclareMathOperator{\Lie}{Lie}
\DeclareMathOperator{\Map}{Map}
\DeclareMathOperator{\cmod}{-mod}
\DeclareMathOperator{\op}{op}
\DeclareMathOperator{\Spec}{Spec}
\newcommand{\cA}{\mathcal{A}}
\newcommand{\cB}{\mathcal{B}}
\newcommand{\cC}{\mathcal{C}}
\newcommand{\cD}{\mathcal{D}}

\newcommand{\cF}{\mathcal{F}}

\newcommand{\cJ}{\mathcal{J}}
\newcommand{\cO}{\mathcal{O}}
\newcommand{\cU}{\mathcal{U}}
\newcommand{\cV}{\mathcal{V}}
\newcommand{\fg}{\mathfrak{g}}
\newcommand{\fh}{\mathfrak{h}}
\newcommand{\fk}{\mathfrak{k}}
\newcommand{\fl}{\mathfrak{l}}
\newcommand{\fq}{\mathfrak{q}}
\newcommand{\bC}{\mathbb{C}}
\newcommand{\bF}{\mathbb{F}}
\newcommand{\bG}{\mathbb{G}}
\newcommand{\bL}{\mathbb{L}}

\newcommand{\bP}{\mathbb{P}}

\newcommand{\bR}{\mathbb{R}}
\newcommand{\bZ}{\mathbb{Z}}

\theoremstyle{plain}
\newtheorem{thm}{Theorem}[subsection]
\newtheorem{cor}[thm]{Corollary}
\newtheorem{lem}[thm]{Lemma}
\newtheorem{prop}[thm]{Proposition}
\theoremstyle{definition}

\newtheorem{cond}[thm]{Condition}
\newtheorem{cons}[thm]{Construction}

\newtheorem{ex}[thm]{Example}
\newtheorem{note}[thm]{Notation}
\newtheorem{rem}[thm]{Remark}
\makeatletter
\newcommand*{\relrelbarsep}{.386ex}
\newcommand*{\relrelbar}{%
	\mathrel{%
		\mathpalette\@relrelbar\relrelbarsep
	}%
}
\newcommand*{\@relrelbar}[2]{%
	\raise#2\hbox to 0pt{$\m@th#1\relbar$\hss}%
	\lower#2\hbox{$\m@th#1\relbar$}%
}
\providecommand*{\rightrightarrowsfill@}{%
	\arrowfill@\relrelbar\relrelbar\rightrightarrows
}

\providecommand*{\xrightrightarrows}[2][]{%
	\ext@arrow 0359\rightrightarrowsfill@{#1}{#2}%
}
\begin{document}
	\title{Dg analogues of the Zuckerman functors and the dual Zuckerman functors I}
	\author{Takuma Hayashi\thanks{Graduate School of Mathematical Sciences, The University of Tokyo, 3-8-1~Komaba, Meguro-ku, Tokyo 153-8914, Japan, htakuma@ms.u-tokyo.ac.jp}}
	\date{}
	\maketitle
	\begin{abstract}
		We study the category of dg Harish-Chandra modules over arbitrary commutative rings, and generalize the induction functor, the production functor, the Zuckerman functor and the dual Zuckerman functor.
	\end{abstract}
	\section{Introduction}
	The theory of Harish-Chandra modules over the field $\bC$ of complex numbers has played a central role in representation theory of real reductive Lie groups as an algebraic model of representations on Hilbert spaces. An important part of this theory is to induce representations cohomologically by the derived functors $\bR I^{\fg,K}_{\fq,M}$ and $\bL P^{\fg,K}_{\fq,M}$. In particular, one can in this way obtain the so-called $A_\fq(\lambda)$-modules. For details, see \cite{MR1330919}.
	
	In this paper, we provide the notion of a differential graded (dg for short) Harish-Chandra module over an arbitrary commutative ring, and generalize the functors $I^{\fg,K}_{\fq,M}$ and $P^{\fg,K}_{\fq,M}$ to this setting.
	\subsection{Motivation}
	The notion of Harish-Chandra modules has been extended in several directions.
	
	One direction is to consider their dg analogues for applications to equivariant derived categories which were introduced by A.\ Beilinson and V.\ Ginzburg (resp.\ J.\ Bernstein and V.\ Lunts) in algebraic (resp.\ geometric) setting in the late 1980s and the 1990s to improve functoriality of equivariant D-modules and sheaves (\cite{MR1237825}, \cite{MR925070}, \cite{MR1317229}). On this course, they also introduced weak $(\fg,K)$-modules and equivariant $(\fg,K)$-complexes. For a uniform approach to both $(\fg,K)$-modules and equivariant complexes of $(\fg,K)$-modules, J.\ Bernstein and V.\ Lunts introduced dg analogues of Harish-Chandra pairs and $(\fg,K)$-modules in \cite{MR1317229}. In that paper, they also introduced the induction functor and the production functor (called the coinduction functor in \cite{MR1317229}). In other words, they constructed left and right adjoint functors of the forgetful functors when the corresponding map of algebraic groups in pairs is the identity map. In \cite{MR2692966}, P.\ Pand\v{z}i\'c explicitly constructed the differential graded analogue of the Zuckerman functor, which he called the equivariant Zuckerman functor, in the setting of equivariant complexes as a right adjoint functor to the forgetful functor.
	
	The second direction is to work over commutative rings and possibly schemes, which has been developed over last ten years and has applications to number theory and mathematical physics.
	
	For example, to get applications to number theory, specifically to rationality and integrality of special values of $L$-functions, one should study rational and integral models of $(\fg,K)$-modules over ground rings like number fields and their rings of integers. In fact, M.\ Harris suggested to work on Harish-Chandra modules and D-modules over number fields to construct models of discrete series representations in \cite{MR3053412} and \cite{10.1093/imrn/rny043}. After this work, G.\ Harder suggested to work over the ring $\bZ$ of integers. He introduced integral models of certain Harish-Chandra modules to estimate contribution of rationality of the Rankin-Selberg $L$-function at the infinite place in \cite{1407.0574}. Moreover, he refined the periods of \cite{Harder2011349} to formulate integrality of special $L$-values. In \cite{MR3770183}, F.\ Januszewski established the cohomological induction over number fields in a fairly similar way to the complex setting. Using his models, he proved rationality of special $L$-values in a representation theoretic way (\cite{1604.04253}, \cite{MR3937337}).
	
	In mathematical physics, J.\ Bernstein et al.\ introduced the contraction families which are Harish-Chandra pairs over the complex projective line $\bP^1$ (\cite{10.1093/imrn/rny146}, \cite{10.1093/imrn/rny147}). In \cite{MR3827131}, E.\ Subag used Harish-Chandra modules over contraction families to reveal the hidden symmetries of the Schr\"odinger equation of the hydrogen atom in two dimensions.
	
	A more primitive example can be obtained from the Lie group SU(1,1). In fact, observe that the Lie algebra
	\[\mathfrak{sl}_2=\bC\left(\begin{array}{cc}
	0&1\\
	0&0\\
	\end{array}
	\right)
	\oplus\bC\left(\begin{array}{cc}
	0&0\\
	-1&0\\
	\end{array}
	\right)\oplus\bC\left(\begin{array}{cc}
	1&0\\
	0&-1\\
	\end{array}
	\right)
	\]
	and the torus
	\[\bG_m=\Spec\bC\left[z^{\pm 1}\right]\cong\left\{\left(\begin{array}{cc}
	z&0\\
	0&z^{-1}\\
	\end{array}
	\right):z\in\bC^\times\right\}\]
	naturally form a Harish-Chandra pair. Replace $\bC$ by $\bZ$ to get a Harish-Chandra pair $(\mathfrak{sl}_2,\Spec\bZ\left[z^{\pm 1}\right])$ over $\bZ$. More precisely, the Lie algebra $\mathfrak{sl}_2$ is spanned by the matrices $\left(\begin{array}{cc}
	0&1\\
	0&0\\
	\end{array}
	\right)$, $\left(\begin{array}{cc}
	0&0\\
	-1&0\\
	\end{array}
	\right)$, and $\left(\begin{array}{cc}
	1&0\\
	0&-1\\
	\end{array}
	\right)$, whose Lie bracket is defined by the usual commutator as matrices. The multiplicative group $\Spec\bZ\left[z^{\pm 1}\right]$ acts on $\mathfrak{sl}_2$ by the conjugation of matrices. The Lie algebra of the multiplicative group $\Spec\bZ\left[z^{\pm 1}\right]$ is naturally isomorphic to the abelian Lie algebra $\bZ$. The $\Spec\bZ\left[z^{\pm 1}\right]$-equivariant Lie algebra homomorphism $\bZ\to\mathfrak{sl}_2$ is given by
	\[1\mapsto \left(\begin{array}{cc}
	1&0\\
	0&-1\\
	\end{array}
	\right).\]
	Such Harish-Chandra pairs are studied in \cite{1712.07336} as split integral models. The (limit of) discrete series representations also admit integral forms, equipped with actions of this integral model of the Harish-Chandra pair (\cite{1712.07336}).
	
	The main purpose of this paper is to generalize the functors $I^{\fg,K}_{\fq,M}$ and $P^{\fg,K}_{\fq,M}$ in these settings at the level of abelian categories in a uniform fashion. They will be tools to construct Harish-Chandra modules over these variant Harish-Chandra pairs. In view of the theory over $\bC$, the next issue is to construct their derived functors to define an analogue of the cohomological induction. In the reductive setting over fields of characteristic 0 without differentials, they can be computed by the Koszul complex (\cite{MR1330919}, \cite{1604.04253}). In the differential graded setting over $\bC$, P.\ Pand\v{z}i\'c constructed an explicit K-projective resolution in \cite{MR2692966} 5.6.5 and \cite{MR2171292} Theorem 3.1 to describe the derived functors. However, these methods do not work in our general setting. In \cite{MR2692966}, he also discussed an abstract homological approach to get resolutions without using the complete reducibility. In \cite{1606.04320}, the theory of model categories is adopted to get the desired derived functors without complicated homological arguments. In fact, we put a model structure, the so called injective model structure, on the category of $(\cA,K)$-modules. We also put another model structure, the so called projective model structure on the same category in case the base ring is a field of characteristic 0 and $K$ is reductive. To relate the approaches from the theories of triangulated categories and model categories, we study the underlying $\infty$-categories in the sense of \cite{luriehigheralgebra}. In the author's paper in preparation, it will be shown in a standard way that they are stable, which means that our $\infty$-category is a homotopical enhancement of the structure of the usual derived categories. The advantage of the $\infty$-categories is that we can use categorical techniques like limits/colimits, the adjoint functor theorem, and generators for our $\infty$-categories which are well developed by J.\ Lurie in \cite{MR2522659} and \cite{luriehigheralgebra}.
	
	A fundamental problem in studying the functor $I^{\fg,K}_{\fq,M}$ is to describe representations obtained by the functor $(\bR)I^{\fg,K}_{\fq,M}$ over commutative rings. The first step in solving this problem is showing compatibility with base change functors (the base change theorem). For example, one may expect that our cohomological induction over subrings of $\bC$ produces forms of $A_\fq(\lambda)$-modules. In view of \cite{MR3770183} and \cite{MR3937337}, they should naively have applications to integrality of special values of automorphic and Rankin-Selberg $L$-functions, though we will need much deeper arguments to prove this. In \cite{MR3853058}, we establish the base change for several practically important situations. On the other hand, we gave nontrivial counterexamples, where the functor $I^{\fg,K}_{\fq,M}$ does not commute with the base change along $\bZ\to\bC$ in \cite{1712.07336}. In fact, we found explicit integral forms of discrete series representations of SU(1,1) which enjoy the principal series type universal property. This is a completely new phenomenon since it never happens over fields of characteristic 0.
	
	The second stage of describing the functor $I^{\fg,K}_{\fg,M}$ will vary upon the result of the first stage and one's interests. In the settings where the base change holds, it is basic to characterize the integral forms in a rigorous way. Suppose that we are given a pair $(\fg,K)$ over an integral domain $k$ with $K$ smooth integral over $k$, a homomorphism $M\to K$ between smooth affine group schemes over $k$, and a $(\fg,M)$-module $V$ which is torsion free over $k$. If the fractional field of $k$ is of characteristic 0 then one can show that $I^{\fg,K}_{\fg,M}(V)$ exhibits the maximal $(\fg,M)$-module whose structure extends to a $(\fg,K)$-module. A more important thing is to investigate the structure of the integral forms or their cohomology in more detail for applications to special $L$-values. For instance, their explicit descriptions are interesting. For deeper applications to number theory, estimates of their rings of definition will be important (see \cite{MR3770183}). In another insight, one of the new things which occur when by considering Harish-Chandra modules over $\bZ$ is torsion, which appears in the derived functor modules over $\bZ$. For example, we formulated and discussed the Borel-Weil-Bott induction over $\bZ$ in \cite{MR3853058}. For a split reductive group $G$ over $\bZ$, it is proved that the algebraic Borel-Weil type induction provides the ``maximal'' $\bZ$-form of the irreducible representations of the complexified group $G\otimes_\bZ\bC$. According to the flat base change theorem, the Borel-Weil-Bott induction over $\bZ$ exhibits the integral model of the Borel-Weil-Bott induction over $\bC$ modulo torsion. In that paper, we found cases when the cohomology has infinite torsion.
	
	In a nontrivial case with the base change failing, a fundamental problem is to describe the resulting modules and new higher derived functor modules as explicitly as possible. It will be interesting to discover the role of these phenomena in representation theory and other branches of mathematics.
	
	In another direction, the theory of D-modules over general base schemes is suggested to be a tool to construct Harish-Chandra modules in \cite{MR3053412}, \cite{10.1093/imrn/rny043}, and \cite{1808.10709}. With the equivariant derived categories of D-modules over number fields (commutative rings), one can obtain the localization of the equivariant Zuckerman functor as in \cite{MR1486143} and \cite{MR2945222}.
	\subsection{Contents of this paper}
	We replace dg Lie algebras $\fg$ by dg associative algebras $\cA$, for example $\cA=U(\fg)$, and we work with Harish-Chandra pairs $(\cA,K)$ (pairs for short) and $(\cA,K)$-modules over an arbitrary commutative ring $k$. Our definitions are straightforward generalizations of those in \cite{MR1330919} and \cite{MR2276617}.
	\begin{note}\label{1.2.1}
		For a pair $(\cA,K)$, write $(\cA,K)\cmod$ for the category of $(\cA,K)$-modules.
	\end{note}
	For a given map $(\cA,K)\to(\cB,L)$ of pairs, one can define the forgetful functor
	\[\cF_{\cB,L}^{\cA,K}:(\cB,L)\cmod\to(\cA,K)\cmod.\]
	\begin{thm}\label{1.2.2}
		\begin{enumerate}
			\renewcommand{\labelenumi}{(\arabic{enumi})}
			\item For a morphism $(\cA,K)\to(\cB,L)$ of pairs, the functor $\cF_{\cB,L}^{\cA,K}$ admits a right adjoint functor $I^{\cB,L}_{\cA,K}$.
			\item For a morphism $(\cA,K)\to(\cB,K)$ of pairs, whose corresponding endomorphism of $K$ is the identity map, the functor $\cF_{\cB,K}^{\cA,K}$ admits a left adjoint functor $P^{\cB,K}_{\cA,K}$.
		\end{enumerate}
	\end{thm}
	This is deduced from its analogue for weak pairs and weak $(\cA,K)$-modules. We introduce the weak concepts in this paper because we can apply generalities on symmetric monoidal categories. More precisely, we will use the fact that cochain complexes of $K$-modules form the closed symmetric monoidal category $K\cmod$. Then a weak pair $(\cA,K)$ is just a monoid object in $K\cmod$, and a weak $(\cA,K)$-module is a left module over the monoid object $(\cA,K)$ in the usual sense in the theory of monoidal categories. Hence their categorical properties and relations easily follow from generalities without any computations. 
	
	Another reason why we adopt this strategy is because of the internal Hom of $K\cmod$ which was used for a construction of the right adjoint functor in \cite{MR1330919}. In the setting of \cite{MR1330919}, the internal Hom of $K\cmod$ is the $K$-finite part of the usual Hom space. However, this description does not make sense in our setting. In fact, the general construction in the proof of \cite{MR2066503} Theorem 1.3.1 is more complicated. Our monoidal category theoretic approach lets us avoid using this explicit description in our construction of the right adjoint functor.
	
	Motivated by Pand\v{z}i\'c's equivariant Zuckerman functor, we summarize the construction of a differential graded analogue of the dual Zuckerman functor in the case when the base ring is $\bC$, and when algebraic groups are reductive in Section 3. The case of $(\fg,K)$-modules without differentials over fields of characteristic 0 has been already done by F.\ Januszewski (\cite{1604.04253} 1.4.2).
	\begin{thm}\label{1.2.3}
		Let $f:(\cA,K)\to(\cA,L)$ be a morphism of pairs, whose corresponding map $\cA\to\cA$ is the identity map. Then the dual Zuckerman functor $\Pi$ and the pseudoforgetful functor $(\cF^\vee)_{\cA,L}^{\cA,K}$ of \cite{MR1330919} naturally extend to an adjunction
		\[\Pi=P^{\cA,L}_{\cA,K}:(\cA,K)\cmod\rightleftarrows(\cA,L)\cmod:(\cF^\vee)_{\cA,L}^{\cA,K}.\]
	\end{thm}
	Combining $P^{\cA,L}_{\cA,K}$ and $(\cF^\vee)_{\cA,L}^{\cA,K}$ with $P^{\cB,K}_{\cA,K}$ (see Theorem \ref{1.2.2} (1)) and $\cF_{\cB,K}^{\cA,K}$ respectively, we obtain a pair ($P^{\cB,L}_{\cA,K}$, $(\cF^\vee)_{\cB,L}^{\cA,K}$) of adjoint functors for an arbitrary map $(\cA,K)\to(\cB,L)$ of pairs. Finally, we show that the principle of induction-in-stages holds for our functors.
	\begin{rem}\label{1.2.4}
		A version of the Zuckerman functor and the dual Zuckerman functor for Lie superalgebras was introduced in \cite{MR1680015}. Our arguments also reprovide Santos' Zuckerman functor and dual Zuckerman functor. Similarly, one can remove the differential graded structures from this paper.
	\end{rem}
	\begin{rem}\label{1.2.5}
		In \cite{MR2692966}, P.\ Pand\v{z}i\'c introduced the notion of triples $(\cA,K,\cD)$ and $(\cA,K,\cD)$-modules to separate two actions of dg algebras. One can define triples $(\cA,K,\cD)$ and $(\cA,K,\cD)$-modules over commutative rings in the manner of \cite{MR2692966}. For full generality, one may allow $\cA$ and $\cD$ to be dg algebras. From the perspectives of the theory of bimodules, it is better to modify the definition so that $\cD$ acts from the right side (see \cite{MR2692966} 3.7 and also 5.4). It also turns out that we need to refine the signs in the compatibility condition of the Lie algebra of $K$. With these minor technical changes in mind, one can easily prove in a similar way to \cite{MR2692966} that for a triple (resp.\ a weak triple) $(\cA,K,\cD)$, the categories of $(\cA,K,\cD)$-modules and $(\cA\otimes \cD^{\op},K)$-modules are naturally isomorphic, where $\cD^{\op}$ is the opposite dg algebra to $\cD$. Hence one can verify the results in this paper for triples via this isomorphism. 
	\end{rem}
	\subsection{Notation}
	Throughout this paper, we assume that there exists a sufficiently large strongly inaccessible cardinal, and fix a Grothendieck universe $\cU$. We freely omit the terminology ``$\cU$-small'' for rings and modules as usual. Let $\bZ$ denote the ring of integers.
	
	For an affine group scheme $K$ over a commutative ring $k$, its coordinate ring will be denoted by $\cO(K)$.
	
	For an object $X$ of a category, let $\id_X$ denote the identity map of $X$. If we are given a pair $F:\cC\to\cD$ and $G:\cD\to\cC$ of adjoint functors between categories $\cC$ and $\cD$, we will denote the unit (resp.\ the counit) by $u_X:X\to G(F(X))$ (resp.\ $v_Y:F(G(Y))\to Y$) for each $X\in\cC$ (resp.\ $Y\in\cD$). In the following, the subscript $X$ of $u_X$ will be omitted if there is no risk of confusion. We will use similar notation for other algebraic structures. 
	
	We next summarize basic terminology and notation used in the theory of symmetric monoidal categories for convenience of the readers who are not familiar with them. See \cite{MR1712872} VII for a reference. A monoidal category is a category $\cV$, equipped with a functor $-\otimes-:\mathcal{V}\times\mathcal{V}\to\mathcal{V}$, an isomorphism $a:(-\otimes-)\otimes -\stackrel{\sim}{\to}-\otimes(-\otimes -)$ of functors from $\mathcal{V}\times\mathcal{V}\times\mathcal{V}$ to $\mathcal{V}$, which is called the associator, an object $I=I_\cV\in\cV$ called the unit, and isomorphisms $r:-\otimes I\stackrel{\sim}{\to}\id_\mathcal{V}:\mathcal{V}\to\mathcal{V}$, $l:I\otimes -\stackrel{\sim}{\to}\id_\mathcal{V}:\mathcal{V}\to\mathcal{V}$ such that the following diagrams are commutative for all quartets $X,Y,Z,W$ of objects in $\cV$:
	\[\begin{tikzcd}
	&(W\otimes X)\otimes (Y\otimes Z)\arrow[rd, "a_{W,X,Y\otimes Z}"]&\\
	((W\otimes X)\otimes Y)\otimes Z \arrow[ru, "a_{W\otimes X,Y,Z}"]\arrow[d, "a_{W,X,Y}\otimes\id_Z"']
	&&W\otimes (X\otimes (Y\otimes Z))\\
	(W\otimes (X\otimes Y))\otimes Z\arrow[rr, "a_{W,X\otimes Y,Z}"]
	&&W\otimes ((X\otimes Y)\otimes Z))\arrow[u, "\id_W\otimes a_{X,Y,Z}"']
	\end{tikzcd}\]
	\[\begin{tikzcd}
	(X\otimes I)\otimes Y\arrow[rd, "r\otimes \id_Y"']\arrow[rr, "a"]&&X\otimes(I\otimes Y)\arrow[ld, "\id_X\otimes l"]\\
	&X\otimes Y.&\\
	\end{tikzcd}\]
	We will sometimes omit the identity map if necessary to save space; for example, we will denote $a_{W,X,Y}\otimes\id_Z$ by $a_{W,X,Y}$ in the above diagram.
	
	For a monoidal category $\cV$, define $-\otimes^{\op}-$ as the composition of $-\otimes-$ and the switch of the components of $\cV\times\cV$:
	\[X\otimes^{\op} Y:=Y\otimes X.\]
	A symmetric monoidal category is a monoidal category $\cV$, equipped with an isomorphism $C:-\otimes-\cong-\otimes^{\op}-$ which satisfies the following properties:
	\begin{enumerate}\renewcommand{\labelenumi}{(\roman{enumi})}
		\item $C_{Y,X}\circ C_{X,Y}=\id_{X\otimes Y}$ for all $X,Y\in\cV$;
		\item The diagram
		\[\begin{tikzcd}
		&X\otimes (Y\otimes Z)\arrow[rd, "C_{X,Y\otimes Z}"]&\\
		(X\otimes Y)\otimes Z\arrow[ru, "a_{X,Y,Z}"]\arrow[d, "C_{X,Y}\otimes\id_Z"']&&(Y\otimes Z)\otimes X\arrow[d, "a_{Y,Z,X}"]\\
		(Y\otimes X)\otimes Z\arrow[rd, "a_{Y,X,Z}"']&&Y\otimes(Z\otimes X)\\
		&Y\otimes(X\otimes Z)\arrow[ru, "\id_Y\otimes C_{X,Z}"']
		\end{tikzcd}\]
		commutes for all $X,Y,Z\in\cV$.
	\end{enumerate}
	We say that a symmetric monoidal category is closed if for every object $X\in\cV$, the functor $-\otimes X:\cV\to\cV$ admits a right adjoint functor $\Map(X,-)$ which is sometimes called the internal Hom. For a pair $X,Y$ of objects in $\cV$, the counit
	\[\Map(X,Y)\otimes X\to Y\]
	will be denoted by $\ev$. Note that the internal Hom is functorial in both the domain and the target (\cite{MR1712872} IV.7 Theorem 3). One can use adjunction to define the enriched composition map
	\[\circ:\Map(Y,Z)\otimes \Map(X,Y)\to\Map(X,Z)\]
	for objects $X,Y,Z\in\cV$ via
	\[\Map(Y,Z)\otimes\Map(X,Y)\otimes X\xrightarrow{\id_{\Map(Y,Z)}\otimes \ev}\Map(Y,Z)\otimes Y\overset{\ev}{\to} Z.\]
	
	Let $\cV$ be a monoidal category. Then a monoid object is an object $A\in\cV$, equipped with two maps $m_A:A\otimes A\to A$ and $j_A:I\to A$ such that the following diagrams are commutative:
	\[\begin{tikzcd}
	&A\otimes(A\otimes A)\arrow[rd, "\id_A\otimes m_A"]&\\
	(A\otimes A)\otimes A\arrow[ru, "a_{A,A,A}"]\arrow[d, "m_A\otimes \id_A"']&&A\otimes A\arrow[d, "m_A"]\\
	A\otimes A\arrow[rr, "m_A"]&&A\\
	\end{tikzcd}\]
	\[\begin{tikzcd}
	I\otimes A\arrow[r, "j_A\otimes \id_A"]\arrow[rd, "l"']&A\otimes A\arrow[d, "m_A"]&A\otimes I\arrow[l, "\id_A\otimes j_A"']\arrow[ld, "r"]\\
	&A.&
	\end{tikzcd}\] 
	When we want to specify the structure morphisms, we will say that $(A,m_A,j_A)$ is a monoid object. We will use similar notation for other algebraic objects. For a monoid object $A\in\cV$, a left $A$-module in $\cV$ is an object $M\in\cV$, equipped with a map $\pi_M:A\otimes M\to M$ such that the following diagrams are commutative:
	\[\begin{tikzcd}
	&A\otimes (A\otimes M)\arrow[rd, "\id_A\otimes \pi_M"]&\\
	(A\otimes A)\otimes M\arrow[ru, "a_{A,A,M}"]\arrow[d, "m_A\otimes \id_M"']&&A\otimes M\arrow[d, "\pi_M"]\\
	A\otimes M\arrow[rr, "\pi_M"]&&M\\
	\end{tikzcd}\]
	\[\begin{tikzcd}
	I\otimes M\arrow[r, "j_A\otimes \id_M"]\arrow[rd, "l"']&A\otimes M\arrow[d, "\pi_M"]&\\
	&M.&
	\end{tikzcd}\] 
	The category of $A$-modules will be denoted by $A\cmod$ in a general setting (Lemma \ref{2.2.5}, Lemma \ref{2.2.6}, and Lemma \ref{2.2.7}).
	
	For example, consider a commutative ring $k$. Write $k\cmod$ for the category of cochain complexes of $k$-modules. Its Hom $k$-module will be denoted by $\Hom$ or $\Hom_k$. This category is closed symmetric monoidal for the usual tensor product $\otimes=\otimes_k$. Let us summarize the structure of this category in detail. Let $M$ be a cochain complex of $k$-modules (or a graded $k$-module). For a homogeneous element $m$, $\bar{m}$ denotes the homogeneous degree of $m$. Write $d_M$ for the differential of $M$. For cochain complexes $M$ and $N$ of $k$-modules, define a new cochain complex $M\otimes N$ as
	\[(M\otimes N)^i=\bigoplus_{p+q=i}M^p\otimes_kN^q\]
	\[d(m\otimes n)=d_Mm\otimes n+(-1)^{\bar{m}}m\otimes d_N n\]
	for any homogeneous element $m\in M$ and any element $n\in N$. The unit object is the complex $k$ which is concentrated in degree $0$ with $k^0=k$. The internal Hom $\Map(M,N)$ is described as
	\[\begin{split}
	\Map(M,N)^i&:=\{\mathrm{graded\ maps\ from}\ M\ \mathrm{to}\ N\left[i\right]\}\\
	&=\prod_p\Hom(M^p,N^{p+i})
	\end{split}\]
	\[df:=d_N\circ f-(-1)^{\bar{f}}f\circ d_M\]
	for a homogeneous element $f\in\Map(M,N)$, where $N\left[i\right]\in k\cmod$ is defined by
	\[(N\left[i\right])^p=N^{p+i}\]
	\[d_{N\left[i\right]}=(-1)^i d_N\]
	with the same $k$-action. The adjunction structure is given as follows: For a cochain complex map $\varphi:M\otimes N\to L$, we define $\alpha(\varphi)\in\Hom(M,\Map(N,L))$ as
	\[\alpha(\varphi)(m)(n)=\varphi(m\otimes n)\]
	for elements $m\in M$ and $n\in N$. Finally, we put the symmetry structure $C_{MN}:M\otimes N\stackrel{\sim}{\to}N\otimes M$ by
	\[m\otimes n\mapsto (-1)^{\bar{m}\bar{n}}n\otimes m\]
	where $m\in M$ and $n\in N$ are homogeneous elements.
	
	A monoid object of $k\cmod$ is a graded $k$-algebra $\cA=\bigoplus_i\cA^i$, equipped with a differential $d$ of $\cA$ of degree 1, i.e., a $k$-module homomorphism $d:\cA\to\cA$ with the following properties:
	\begin{enumerate}\renewcommand{\labelenumi}{(\roman{enumi})}
		\item $d^2=0$;
		\item $d(\cA^i)\subset \cA^{i+1}$;
		\item $d(ab)=(da)b+(-1)^{\bar{a}}a(db)$ for all homogeneous elements $a,b\in\cA$.
	\end{enumerate}
	This is known as a dg algebra over $k$. We will simply say that $\cA$ is a dg algebra. For a dg algebra $\cA$, a left $\cA$-module in $k\cmod$ is a graded left $\cA$-module $M$, equipped with a differential $d$ of $M$ of degree 1, i.e., a $k$-module homomorphism $d:M\to M$ satisfying the following conditions:
	\begin{enumerate}
		\renewcommand{\labelenumi}{(\roman{enumi})}
		\item $d^2=0$;
		\item $d(M^i)\subset M^{i+1}$;
		\item $d(am)=(da)m+(-1)^{\bar{a}}a(dm)$ for any homogeneous element $a\in\cA$ and any element $m\in M$.
	\end{enumerate}
	This is known as a left dg $\cA$-module. We will simply say that $M$ is a dg (left) $\cA$-module. A morphism of dg $\cA$-modules is defined as a homomorphism of $k\cmod$ respecting the action maps of $\cA$.
	
	For a digression, recall that a dg Lie algebra is a dg $k$-module $\fg$, equipped with a graded bilinear map $\left[-,-\right]:\fg\otimes\fg\to\fg$ which satisfies the following conditions for any three homogeneous elements $x, y, z\in\fg$:
	\begin{enumerate}
		\renewcommand{\labelenumi}{(\roman{enumi})}
		\item $\left[x,y\right]=-(-1)^{\bar{x}\bar{y}}\left[y,x\right]$;
		\item $\left[x,\left[y, z\right]\right] = \left[\left[x, y \right], z\right] + (-1)^{\bar{x}\bar{y}}\left[y,\left[x, z\right]\right]$;
		\item $d\left[x, y\right]=\left[dx,y\right]+(-1)^{\bar{x}}\left[x, dy\right]$.
	\end{enumerate}
	For example, a Lie algebra over $k$ is regarded as a dg Lie algebra concentrated in degree 0.
	
	We next consider an affine group scheme $K$ over $k$. Then the category of cochain complexes of representations of $K$ will be denoted by $K\cmod$. This is closed symmetric monoidal for the usual tensor product $\otimes=\otimes_k$. Its objects are sometimes called dg $K$-modules or dg representations of $K$. Note that if $K$ is flat over $k$, $K\cmod$ is a Grothendieck abelian category. Namely, $K\cmod$ is a locally small cocomplete abelian category with generators, which satisfies the property that filtered colimits are exact. Moreover, its monomorphisms are the (degreewise) injective homomorphisms. See \cite{MR2138086} Proposition 1.2 for example.
	
	If we are given an affine group scheme with a capital symbol, its Lie algebra will be denoted by the corresponding German symbol. Let $K$ and $L$ be affine group schemes over $k$. For a dg $K$-module $M$ with an action map $\nu$, the corresponding differential dg representation will be denoted by $d\nu$. For an affine group scheme homomorphism $f:K\to L$, its differential will be denoted by $df:\fk\to\fl$; it is a Lie algebra homomorphism.
	
	We now return to general monoidal categories. A lax monoidal functor is a functor $F:\cV\to\cV'$, equipped with a morphism $\mu_0:I_{\cV'}\to F(I_\cV)$ and a natural transformation $\mu:F(-)\otimes F(-)\to F(-\otimes -)$ such that the following diagrams are commutative for all $X,Y,Z\in\cV$:
	\[\begin{tikzcd}
	(F(X)\otimes F(Y))\otimes F(Z)\ar[rr, "a"]\ar[d, "\mu_{X,Y}\otimes\id_{F(Z)}"']
	&&F(X)\otimes (F(Y)\otimes F(Z))\ar[d, "\id_{F(X)}\otimes\mu_{Y,Z}"]\\
	F(X\otimes Y)\otimes F(Z)\ar[d, "\mu_{X\otimes Y,Z}"']&& F(X)\otimes F(Y\otimes Z)\ar[d, "\mu_{X,Y\otimes Z}"]\\
	F((X\otimes Y)\otimes Z)\ar[rr, "F(a_{X,Y,Z})"]&&F(X\otimes (Y\otimes Z))
	\end{tikzcd}\]
	\[\begin{tikzcd}
	I_{\cV'}\otimes F(X)\ar[r, "l"]\ar[d, "\mu_0\otimes\id_{F(X)}"']&F(X)\\
	F(I_\cV)\otimes F(X)\ar[r, "\mu"]&F(I_\cV\otimes X)\ar[u, "F(l)"']
	\end{tikzcd}\]
	\[\begin{tikzcd}
	F(X)\otimes I_{\cV'}\ar[r, "r"]\ar[d, "\id_{F(X)}\otimes\mu_0"']&F(X)\\
	F(X)\otimes F(I_\cV)\ar[r, "\mu"]&F(X\otimes I_\cV).\ar[u, "F(r)"']
	\end{tikzcd}\]
	We say that $F$ is monoidal if the map $\mu_0:I_{\cV'}\to F(I_\cV)$ and the natural transformation $\mu:F(-)\otimes F(-)\to F(-\otimes -)$ are isomorphisms. If $\cV$ and $\cV'$ are symmetric, a monoidal functor $F:\cV\to\cV'$ is said to be symmetric if the diagram
	\[\begin{tikzcd}
	F(X)\otimes F(Y)\ar[rr, "C_{F(X),F(Y)}"]\ar[d, "\mu"']&&F(Y)\otimes F(X)\ar[d, "\mu"]\\
	F(X\otimes Y)\ar[rr, "F(C_{X,Y})"']&&F(Y\otimes X)
	\end{tikzcd}\]
	commutes for all $X,Y\in\cV$.
	\renewcommand{\abstractname}{Acknowledgments}
	\begin{abstract}
		First of all, I would like to express my gratitude to my advisor Professor Hisayosi Matumoto. He always listens to my talks with much time, and he gives me much useful advice. Especially, he introduced the papers \cite{MR2692966}, \cite{MR2171292} and \cite{MR2276617} to me. I am grateful to Professor Pavle Pand\v{z}i\'c for his work on (weak) $(\cA,K,\cD)$-modules and equivariant Zuckerman functor in his thesis \cite{MR2692966}. This was a starting point for my study in the present paper. I am also indebted to him for his polite answers to my questions during a coffee break in the conference ``Representations of reductive groups: A conference dedicated to David Vogan on his 60th birthday'' at MIT. I would like to thank the referee for suggesting an idea to simplify Section 3. I would like to thank Tobias Columbus for telling me about \cite{MR2066503} Theorem 1.3.1. I also thank Tobias Columbus and Jun Yoshida for discussions on the closed structure on the symmetric monoidal category of comodules over a flat commutative Hopf algebra. Finally, I would like to thank all my colleagues, especially, Tomoki Mihara, Yoshiki Oshima, Hironori Oya, Koji Shimizu, Yuichiro Tanaka, and Kohei Yahiro for teaching me various parts of mathematics and discussing many topics with me.
		
		This work was supported by JSPS Kakenhi Grant Number JP15J06457.
	\end{abstract}
	\section{Basic theory of dg Harish-Chandra modules}
	\subsection{($\cA,K$)-modules and their weak analogues}
	A weak pair $(\cA,K)$ is a pair of a flat affine group scheme $K$ over $k$ and a monoid object $\cA$ in the symmetric monoidal category $K\cmod$. In other words, a weak pair consists of a flat affine group scheme $K$ and a dg $k$-algebra $(\cA,m_\cA,j_\cA)$, equipped with a $K$-action $\phi$, namely, a dg $K$-module structure on $\cA$ with the properties that $m_\cA$ and $j_\cA$ are $K$-equivariant. For weak pairs $(\cA,K)$ and $(\cB,L)$, a weak map from $(\cA,K)$ to $(\cB,L)$ is a pair $f=(f_a,f_k)$ consisting of a dg algebra homomorphism $f_a:\cA\to\cB$ and a group scheme homomorphism $f_k:K\to L$ with the property that $f_a$ is $K$-equivariant via $f_k$.
	
	To introduce the notion of pairs, recall how the adjoint ``representation'' is constructed: For any commutative $k$-algebra $R$, we have a homomorphism
	\[R\left[\epsilon\right]/(\epsilon^2)\to R;\ a+b\epsilon\mapsto a.\]
	Using this map, we define a functor $\Lie K$ from the category of commutative $k$-algebras to that of groups as 
	\[(\Lie K)(R)=\Ker(K(R\left[\epsilon\right]/(\epsilon^2))\to K(R)).\]
	We now define an action Ad of $K(R)$ on $(\Lie K)(R)$ by conjugation via the map
	\[R\to R\left[\epsilon\right]/(\epsilon^2);\ a\mapsto a+0\epsilon\]
	for each $R$. Hence we obtain a homomorphism of group functors
	\[\Ad:K\to\Aut(\Lie K).\]
	One can obtain the natural structure of a Lie algebra over $R$ on $(\Lie K)(R)$ from this action (see \cite{MR563524} Chapter II, \S 4, 4.2 and 4.5 Proposition). In particular, the Lie algebra $\fk$ of $K$ is defined by $\fk=(\Lie K)(k)$.
	
	To obtain a representation of $K$ on the $k$-module $\fk$ from this construction, consider the following condition for an affine group scheme $K$:
	\begin{cond}\label{2.1.1}
		The conormal $k$-module $I_e/I_e^2$ is finitely generated and projective, where $I_e$ denotes the kernel of the counit map of the coordinate ring of $K$.
	\end{cond}
	\begin{ex}[\cite{MR563524} Chapter II, \S 4, 4.8]\label{2.1.2}
		Condition \ref{2.1.1} is satisfied in either of the following cases:
		\begin{enumerate}
			\renewcommand{\labelenumi}{(\roman{enumi})}
			\item $K$ is smooth over $k$;
			\item $k$ is a field, and $K$ is of finite type over $k$.
		\end{enumerate}
	\end{ex}
	\begin{lem}[\cite{MR563524} Chapter II, \S 4, 4.8 Proposition]\label{2.1.3}
		Let $K$ be an affine group scheme over $k$. Then the $R$-module homomorphism
		\[R\otimes_k(\Lie K)(k)\to(\Lie K)(R)\]
		is an isomorphism for every commutative $k$-algebra $R$ if and only if $K$ satisfies Condition \ref{2.1.1}.
	\end{lem}
	\begin{ex}\label{2.1.4}
		The additive group $\bZ/2\bZ$ attaches the diagonalizable group $K=\Spec\bZ\left[t\right]/(t^2-1)$. Then we have $(\Lie K)(\bZ)=0$ and $(\Lie K)(\bF_2)=\bF_2$, where $\bF_2$ is the finite field of two elements. In \cite{MR3853058}, we deal with the case when $K$ has not so bad singularity, for example, $K$ is flat and finitely presented over $k$.
	\end{ex}
	Therefore if $K$ satisfies Condition \ref{2.1.1}, the adjoint representation $\Ad$ of $K$ on the $k$-module $\fk=(\Lie K)(k)$ is obtained by
	\[K(R)\times (R\otimes_k (\Lie K)(k))\to K(R)\times (\Lie K)(R)\overset{\Ad}{\to}(\Lie K)(R)\cong R\otimes_k(\Lie K)(k),\]
	where $R$ runs through all commutative $k$-algebras. We are now ready to define pairs: Let $(\cA,K,\phi)$ be a weak pair with the property that $K$ satisfies Condition \ref{2.1.1}, and let $\psi:\fk\to\cA$ be a $K$-equivariant dg Lie algebra homomorphism. Then $(\cA,K)$ is called a pair if the following equality holds for any $\xi\in\fk$:
	\[d\phi(\xi)=\left[\psi(\xi),-\right]:\cA\to\cA.\]
	A map
	\[(\cA,K,\phi_\cA,\psi_\cA)\to(\cB,L,\phi_\cB,\psi_\cB)\]
	of pairs is a weak map $f=(f_a,f_k)$ respecting $\psi$, i.e., satisfying the equality
	\[f_a\circ\psi_\cA=\psi_\cB\circ df_k:\fk\to\cB.\]
	\begin{ex}\label{2.1.5}
		Let $K$ be a flat affine group scheme.
		\begin{enumerate}
			\renewcommand{\labelenumi}{(\arabic{enumi})}
			\item We put $\cA=k$ with the trivial action. This naturally gives a weak pair $(k,K)$. We call it the trivial weak pair. The trivial Lie algebra homomorphism $\fk\to k$ defines a structure of a pair on $(k,K)$ if $K$ satisfies Condition \ref{2.1.1}.
			\item Let $(\cA,K)$ be a weak pair. Then we have a unique weak map $f=(f_a,\id_K)$ from the trivial weak pair $(k,K)$ to $(\cA,K)$.
			\item Suppose that $K$ satisfies Condition \ref{2.1.1}. We set $\cA$ as the enveloping algebra $U(\fk)$ with the adjoint action. Then $(U(\fk),K)$ forms a pair by the canonical Lie algebra homomorphism $\fk\to U(\fk)$. We call it the trivial pair.
			\item Let $(\cA,K)$ be a pair. Then we have a unique map $f=(f_a,\id_K)$ from the trivial pair $(U(\fk),K)$ to $(\cA,K)$.
			\item For a weak pair $(\cA,K)$, $(\cA^{\op},K)$ is a weak pair where $\cA^{\op}$ denotes the opposite dg algebra to $\cA$ with the same $K$-action. We call it the opposite weak pair to $(\cA,K)$. This is the opposite monoid to $(\cA,K)$ in the sense of monoidal category theory.
			\item For a pair $(\cA,K,\phi,\psi)$, $(\cA^{\op},K,\phi,-\psi)$ is a pair. We call it the opposite pair to $(\cA,K,\phi,\psi)$.
			\item For two weak pairs $(\cA,K)$ and $(\cB,K)$, $(\cA\otimes\cB,K)$ naturally forms a weak pair. This follows from general theory of monoids in symmetric monoidal categories.
			\item For two pairs $(\cA,K,\phi_\cA,\psi_\cA)$ and $(\cB,K,\phi_\cB,\psi_\cB)$,
			\[(\cA\otimes\cB,K,\phi_\cA\otimes\phi_\cB,\psi_\cA\otimes 1+1\otimes\psi_\cB)\]
			is a pair.
		\end{enumerate}
	\end{ex}
	For a weak pair $(\cA,K)$, the category of (left) weak $(\cA,K)$-modules is defined as the category of left modules over $(\cA,K)$ in the sense of 1.3. In other words, a weak $(\cA,K)$-module is a dg $k$-module $M$ with a left dg $\cA$-module structure $\pi$ and a $K$-module structure $\nu$ satisfying the condition that $\pi$ is $K$-equivariant, i.e., the following diagram commutes:
	\[\begin{tikzcd}
	K\times\cA\otimes M\ar[rr, "\id_K\times\pi"]\ar[d, "\phi\otimes\nu"']\ar[d]
	&&K\times M\ar[d, "\nu"]\\
	\cA\otimes M\ar[rr, "\pi"']&&M.
	\end{tikzcd}\]
	For a pair $(\cA,K)$, we define an $(\cA,K)$-module as a weak $(\cA,K)$-module $M$ satisfying $\pi(\psi(\xi))=d\nu(\xi)$ for every $\xi\in\fk$.
	\begin{note}\label{2.1.6}
		\begin{enumerate}
			\renewcommand{\labelenumi}{(\arabic{enumi})}
			\item For a weak pair $(\cA,K)$, write $(\cA,K)\cmod_w$ for the category of weak $(\cA,K)$-modules.
			\item For a pair $(\cA,K)$, write $(\cA,K)\cmod$ for the full subcategory of $(\cA,K)\cmod_w$ consisting of $(\cA,K)$-modules.
		\end{enumerate}
	\end{note}
	\begin{lem}\label{2.1.7}
		Let $K$ be a flat affine group scheme. Then the category of weak $(k,K)$-modules is canonically isomorphic to $K\cmod$.
	\end{lem}
	\begin{proof}
		Observe that the trivial $K$-module $k$ is the unit object of the monoidal category $K\cmod$. The assertion is now obvious.
	\end{proof}
	\subsection{Functors}
	Our first goal is to give dg analogues of functors $\cF_{\fg,K}^{\fh,L},P^{\fg,K}_{\fh,K}$ and $I^{\fg,K}_{\fh,L}$ (\cite{MR1330919}). The first proposition below is obvious.
	\begin{prop}\label{2.2.1}
		Let $f=(f_a,f_k):(\cA,K)\to(\cB,L)$ be a weak map of weak pairs. For a weak $(\cB,L)$-module $(M,\pi_2,\nu_2)$, the dg $k$-module $M$ admits a weak $(\cA,K)$-module structure $(\pi_1,\nu_1)$ as follows:
		\[\pi_1=\pi_2\circ f_a\]
		\[\nu_1=\nu_2\circ f_k.\]
		Moreover, if $f$ is a map of pairs and $M$ is a $(\cB,L)$-module, the resulting weak $(\cA,K)$-module is actually an $(\cA,K)$-module.
	\end{prop}
	As a consequence, we obtain the following two forgetful functors:
	\[\cF_{\cB,L,w}^{\cA,K}:(\cB,L)\cmod_w\to(\cA,K)\cmod_w;\]
	\[\cF_{\cB,L}^{\cA,K}:(\cB,L)\cmod\to(\cA,K)\cmod.\]
	We construct their right adjoint functors. For a pair $(\cA,K)$, let
	\[\cJ_{\cA,K}:(\cA,K)\cmod\to(\cA,K)\cmod_w\]
	denote the natural fully faithful embedding. The next lemma reduces existence of the right adjoint functor $I^{\fg,K}_{\fq,L}$ to the cases of weak modules.
	\begin{lem}[\cite{MR2692966} 5.7.3]\label{2.2.2}
		Let $(\cA,K)$ be a pair. Then the full subcategory $(\cA,K)\cmod$ of $(\cA,K)\cmod_w$ is both a localization and a colocalization, i.e., $\cJ_{\cA,K}$ admits both a left adjoint functor $(-)_\fk$ and a right adjoint functor $(-)^\fk$. Moreover, the adjoint functors only depend on $K$ in the sense that we have the following isomorphisms for a map $(\cA,K)\to(\cB,K)$ of pairs, whose corresponding map $K\to K$ is the identity map:
		\[(-)^\fk\circ\cF^{\cA,K}_{\cB,K,w}\cong\cF^{\cA,K}_{\cB,K,w}\circ(-)^\fk;\]
		\[(-)_\fk\circ\cF^{\cA,K}_{\cB,K,w}\cong\cF^{\cA,K}_{\cB,K,w}\circ(-)_\fk.\]
	\end{lem}
	\begin{proof}
		Let $(\cA,K,\phi,\psi)$ be a pair, and $(M,\pi,\nu)$ be a weak $(\cA,K)$-module. Set
		\[\omega(\xi)=d\nu(\xi)-\pi(\psi(\xi))\]
		for $\xi\in\fk$. Then $\omega$ is a dg representation of $\fk$ by \cite{MR1486143} 2.1. Moreover, $\omega$ commutes with the action of $\cA$ and it is equivariant with respect to the action of $K$. Hence the dg submodule
		\[N:=\omega(\fk)M\]
		inherits a weak $(\cA,K)$-submodule structure, and $N$ is stable under the action $\omega$. Hence the quotient dg module $M_\fk=M/N$ defines a left adjoint functor of $\cJ_{\cA,K}$. Similarly, let $M^\fk\subset M$ be the submodule where $\omega$ acts trivially. Passing to the adjunction, one can identify it with the kernel of the $K$-module homomorphism $M\to\Map(\fk,M)$ (\cite{MR2015057} I.2.7 (5)). Hence $M^\fk$ is an $(\cA,K)$-submodule, and $(-)^\fk$ gives rise to a right adjoint functor to $\cJ_{\cA,K}$.
	\end{proof}
	In view of the classical construction of the functor $I^{\fg,K}_{\fq,L}$ (especially in the case $K=L$), we should use the internal Hom of $K\cmod$ and $L\cmod$. Let us recall the construction:
	\begin{cons}\label{2.2.3}
		Let $C$ be a (flat) coalgebra over $k$, $V$ be a right dg $C$-comodule, $\rho=\rho_V:V\to V\otimes C$ be its coaction, and $\bar{V}$ be the cokernel of $\rho$. Since $\rho$ can be regarded as a dg $C$-comodule homomorphism from the given right dg $C$-comodule $V$ to the free right dg $C$-comodule $V\otimes C$, $\bar{V}$ is canonically equipped with the structure of a right dg $C$-comodule, whose coaction will be denoted by $\bar{\rho}$. We again regard $\bar{\rho}$ as a dg $C$-comodule homomorphism from $\bar{V}$ to the free right dg $C$-comodule $\bar{V}\otimes C$. Compose $\bar{\rho}$ with the canonical quotient map $\coker\rho:V\otimes C\to \bar{V}$ to get a $C$-homomorphism $\tau_V:V\otimes C\to\bar{V}\otimes C$ between free right dg $C$-comodules, whose kernel (as a $k$-module) is $V$ since the coaction maps are injective.
	\end{cons}
	Recall that for a flat affine group scheme $K$, $K\cmod$ can be identified with the symmetric monoidal category of right dg comodules over the coordinate ring $H:=\cO(K)$ of $K$ (\cite{MR547117} Theorem 3.2). In the following, we will denote the Hom bifunctor of right dg $H$-comodules by $\Hom_H$. Let $m$ (resp.\ $\epsilon$) be the multiplication map (resp.\ counit) of $H$.
	\begin{prop}[\cite{MR2066503} Theorem 1.3.1]\label{2.2.4}
		Let $M,M'$ be right dg comodules over $H$, and $N,N'$ be dg $k$-modules.
		\begin{enumerate}
			\renewcommand{\labelenumi}{(\arabic{enumi})}
			\item Regard $\Map(M,N)\otimes H$ as a free right dg $C$-comodule. Then there is a bijection $\gamma:\Hom_H(M',\Map(M,N)\otimes H)\cong\Hom_H(M'\otimes M,N\otimes H)$ which is natural in $M$ and $M'$.
			\item For a comodule homomorphism $f:N\otimes H\to N'\otimes H$, set $\hat{f}_M=\gamma^{-1}(f\circ \gamma(\id_{\Map(M,N)\otimes H})):\Map(M,N)\otimes H\to\Map(M,N')\otimes H$. Then the diagram
			\[\begin{tikzcd}
			\Hom_H(M',\Map(M,N)\otimes H)\ar[r, "\hat{f}_M\circ-"]\ar[d, "\gamma"']&\Hom_H(M',\Map(M,N')\otimes H)\ar[d, "\gamma"]\\
			\Hom_H(M'\otimes M,N\otimes H)\ar[r, "f\circ-"]&\Hom_H(M'\otimes M,N'\otimes H)
			\end{tikzcd}\]
			commutes.
			\item The homomorphism $\gamma(\id_{\Map(M,N)\otimes H})$ is given by the composite map
			\[\begin{split}
			\Map(M,N)\otimes H\otimes M
			&\xrightarrow{\rho_M}\Map(M,N)\otimes H\otimes M\otimes H\\
			&\xrightarrow{C_{H,M}}\Hom_k(M,N)\otimes M\otimes H\otimes H\\
			&\xrightarrow{\ev\otimes m} N\otimes H.
			\end{split}
			\]
			\item The internal Hom of $K\cmod$ is given by
			\[F(M,M')=\Ker (\hat{(\tau_{M'})}_M:\Map(M,M')\otimes H \to \Map(M,\bar{M}')\otimes H).\]
		\end{enumerate}
	\end{prop}
	\begin{proof}
		Part (1) is due to the following formal computation of the usual adjunctions:
		\[\begin{split}
		\Hom_H(M',\Map(M,N)\otimes H)
		&\cong\Hom(M',\Map(M,N))\\
		&\cong\Hom(M'\otimes M,N)\\
		&\cong\Hom_H(M'\otimes M,N\otimes H).
		\end{split}\]
		It is clear that this sequence of bijections is natural in $M$ and $M'$.
		
		Part (2) immediately follows from the naturality of (1).
		
		We unwind the definitions to show (3). Firstly, the identity map $\id_{\Map(M,N)\otimes H}$ goes to
		\[\id_{\Map(M,N)}\otimes\epsilon:\Map(M,N)\otimes H\to\Map(M,N).\]
		Pass to the adjunction of $\Map$ and $\otimes$ to get a map
		\[\Map(M,N)\otimes H\otimes M\to N;\ g\otimes a\otimes m\mapsto \epsilon(a)g(m).\]
		We next regard $\Map(M,N)\otimes H\otimes M$ as the tensor comodule of $\Map(M,N)\otimes H$ and $M$. The map $\gamma(\id_{\Map(M,N)\otimes H})$ is obtained by universally extending to a comodule homomorphism $\Map(M,N)\otimes H\otimes M\to N\otimes H$. Explicitly, $\gamma(\id_{\Map(M,N)\otimes H})$ is expressed by
		\[\begin{split}
		g\otimes a\otimes m
		&\mapsto \sum g\otimes a_1\otimes m_1\otimes a_2c_2\\
		&\mapsto \sum \epsilon(a_1)g(m_1)\otimes a_2c_2\\
		&=\sum g(m_1)\otimes ac_2
		\end{split}\]
		(the Sweedler notation). This coincides with the map in the assertion.
		
		The general adjunction is obtained from the isomorphism
		\[\begin{split}
		\Hom_H(-\otimes M,M')
		&\cong\Ker\Hom_H(-\otimes M,\tau_{M'})\\
		&\cong\Ker\Hom_H(-,\hat{(\tau_{M'})}_M)\\
		&\cong\Hom_H(-,F(M,M')).
		\end{split}\]
		This shows (4).
	\end{proof}
	We now factorize a given morphism $f=(f_a,f_k):(\cA,K)\to (\cB,L)$ of weak pairs as
	\[(\cA,K)\xrightarrow{(f_a,\id_K)} (\cB,K)\xrightarrow{(\id_\cB,f_k)} (\cB,L)\]
	to reduce construction of $P^{\cB,K}_{\cA,K,w}$ and $I^{\cB,K}_{\cA,L,w}$ to the following three well-known results in the theory of monoidal categories:
	\begin{lem}\label{2.2.5}
		Let $F:\cV\to\cV'$ be a monoidal functor between monoidal categories with a right adjoint functor $G$, and $(A,m)$ be a monoid object of $\cV$. Then the adjunction $(F,G)$ extends to
		\[F:A\cmod\rightleftarrows F(A)\cmod:G.\]
	\end{lem}
	\begin{proof}
		Since $F$ is monoidal, $F(A)$ is a monoid object of $\cV'$ in a natural way. Observe also that $G$ is lax monoidal in a canonical way (\cite{MR2724388} Proposition 3.84). For an $F(A)$-module $(N,\pi_N)$, define a map in $\cV$ as
		\[A\otimes G(N)\xrightarrow{u_A\otimes\id_{G(N)}} G(F(A))\otimes G(N)\xrightarrow{\mu_{F(A),N}} G(F(A)\otimes N)\xrightarrow{G(\pi_N)} G(N).\]
		This determines the structure of a left $A$-module on $G(N)$. In fact, use the lax structure to get the following commutative diagrams:
		\[\begin{tikzcd}
		A\otimes (A\otimes G(N))\ar[r, "u_A\otimes u_A"]\ar[d, "\id_A\otimes u_A\otimes\id_{G(N)}"']
		&G(F(A))\otimes (G(F(A))\otimes G(N))\ar[d, "\id\otimes G(\mu)"]&\\
		A\otimes (G(F(A))\otimes G(N))\ar[d, "\id_A\otimes\mu"']
		&G(F(A))\otimes G(F(A)\otimes N)\ar[lddd, "\id_{G(F(A))}\otimes G(\pi_N)"]\ar[dd, "\mu_{F(A), F(A)\otimes N}"]\\
		A\otimes G(F(A)\otimes N)\ar[d, "G(\pi_N)"']\\
		A\otimes G(N)\ar[d, "u_A\otimes\id_{G(N)}"']&G(F(A)\otimes (F(A)\otimes N))\ar[ldd, "G(\id_{F(A)}\otimes\pi_N)"]\\
		G(F(A))\otimes G(N)\ar[d, "\mu_{F(A),N}"']\\
		G(F(A)\otimes N)\ar[d, "G(\pi_N)"']\\
		G(N)
		\end{tikzcd}\]
		\[\begin{tikzcd}
		(A\otimes A)\otimes G(N)\ar[rr, "(u_A\otimes u_A)\otimes\id_{G(N)}"]\ar[ddd, "m_A\otimes\id_{G(N)}"']
		\ar[rdd, "u_{A\otimes A}\otimes\id_{G(N)}"']
		&&(G(F(A))\otimes G(F(A)))\otimes G(N)\ar[ld, "\mu_{F(A),F(A)}\otimes\id_{G(N)}"]\\
		&G(F(A)\otimes F(A))\otimes G(N)\ar[r, "\mu_{F(A)\otimes F(A),N}"']\ar[d, "G(\mu_{A,A})\otimes\id_{G(N)}"]
		&G((F(A)\otimes F(A))\otimes N)\ar[d, "G(\mu_{A,A}\otimes\id_N)"]\\
		&G(F(A\otimes A))\otimes G(N)\ar[r, "\mu_{F(A\otimes A),N}"']\ar[d, "G(F(m_A))\otimes G(\id_N)"']
		&G(F(A\otimes A)\otimes N)\ar[d, "G(F(m_A)\otimes \id_N)"]\\
		A\otimes G(N)\ar[r, "u_A\otimes\id_{G(N)}"']&G(F(A))\otimes G(N)\ar[r, "\mu_{F(A),N}"']&G(F(A)\otimes N)\ar[d, "G(\pi_N)"]\\
		&&G(N)
		\end{tikzcd}\]
		\[\begin{tikzcd}
		A\otimes G(N)\ar[r, "u_A"]&G(F(A))\otimes G(N)\ar[r, "\mu_{F(A),N}"]&G(F(A)\otimes N)\ar[r, "G(\pi_N)"]&G(N)\\
		I_\cV\otimes G(N)\ar[u, "j_A\otimes\id_{G(N)}"]\ar[rd, "\mu_0"']\ar[r, "u_{I_\cV}"]\ar[d, "l"']
		&G(F(I_\cV))\otimes G(N)\ar[r, "\mu_{F(I_\cV),N}"]\ar[u, "G(F(j_A))\otimes\id_{G(N)}"]\ar[d, "G(\mu_0^{-1})"']
		&G(F(I_\cV)\otimes N)\ar[d, "G(\mu_0^{-1}\otimes\id_N)"']\ar[u, "G(F(j_A)\otimes\id_N)"]&\\
		G(N)\ar[rd, "G(l^{-1})"']&G(I_{\cV'})\otimes G(N)\ar[d, "\mu"]\ar[r, "\mu"']&G(I_{\cV'}\otimes N)\ar[ruu, "G(l)"']\\
		&G(I_{\cV'}\otimes N)\ar[ru, equal]
		\end{tikzcd}\]
		(use the colax-lax property of \cite{MR2724388} Definition 3.81 and Proposition 3.82). Compare the first and second diagrams to get the associativity of the action. The third diagram shows the unitality. This construction is clearly functorial. The proof will be completed by seeing that the counit and the unit are compatible with actions of $A$ and $F(A)$ respectively. It is again shown by routine diagram chases for an $A$-module $M$ and an $F(A)$-module $N$:
		\[\begin{tikzcd}
		A\otimes M\ar[ddd, "\pi_M"']\ar[rr, "\id_A\otimes u_M"]\ar[rdd, "u_{A\otimes M}"']
		&&A\otimes G(F(M))\ar[d, "u_A\otimes\id_{G(F(M))}"]\\
		&&G(F(A))\otimes G(F(M))\ar[d, "\mu_{F(A),F(M)}"]\\
		&G(F(A\otimes M))\ar[rd, "G(F(\pi_M))"']\ar[r, "G(\mu^{-1}_{A,M})"]&G(F(A)\otimes F(M))\ar[d, "G(\pi_{F(M)})"]\\
		M\ar[rr, "u_M"']&&G(F(M))
		\end{tikzcd}\]
		\[\begin{tikzcd}
		F(A)\otimes F(G(N))\ar[d, "\mu_{A,G(N)}"']\ar[rr, "\id_{F(A)}\otimes v_N"]&&F(A)\otimes N\ar[dddd, "\pi_N"]\\
		F(A\otimes G(N))\ar[d, "F(u_A\otimes\id_{G(N)})"']&F(G(F(A)))\otimes F(G(N))\ar[ru, "v_{F(A)}\otimes v_N"]\\
		F(G(F(A))\otimes G(N))\ar[d, "F(\mu_{F(A),N})"']\ar[ru, "\mu_{G(F(A)),G(N)}^{-1}"']\\
		F(G(F(A)\otimes N))\ar[d, "F(G(\pi_N))"']\arrow[rruuu,bend right,"v_{F(A)\otimes N}"']\\
		F(G(N))\ar[rr, "v_N"]&&N.
		\end{tikzcd}\]
	\end{proof}
	\begin{lem}\label{2.2.6}
		Let $\cV$ be a closed symmetric monoidal category with equalizers, and $f:A\to B$ be a morphism of monoids of $\cV$. Then the forgetful functor
		\[B\cmod\to A\cmod\]
		admits a right adjoint functor.
	\end{lem}
	\begin{proof}
		We ignore the associator $a$ in this proof for convenience.
		
		Define $\hat{m}^{\op}:B\to\Map(B,B)$ by the opposite multiplication
		\[B\otimes B\overset{C}{\cong}B\otimes B\overset{m}{\to}B.\]
		
		Let $(M,\pi_M)$ be an $A$-module. Then $\Map(B,M)$ is a left $B$-module for the multiplication on the domain $B$ from the right side:
		\[\begin{split}
		B\otimes \Map(B,M)&\overset{\hat{m}^{\op}}{\to}\Map(B,B)\otimes\Map(B,M)\\
		&\overset{C}{\cong}\Map(B,M)\otimes\Map(B,B)\\
		&\overset{\circ}{\to}\Map(B,M).
		\end{split}\]
		This corresponds to the composite map
		\[\begin{split}
		B\otimes\Map(B,M)\otimes B
		&\xrightarrow{C_{B,\Map(B,M)\otimes B}}\Map(B,M)\otimes B\otimes B\\
		&\overset{m}{\to}\Map(B,M)\otimes B\\
		&\overset{\ev}{\to} M.
		\end{split}\]
		To see the associativity, observe that the following diagram commutes:
		\[\begin{tikzcd}
		B\otimes B\ar[r, "m"]\ar[d, "\hat{m}^{\op}\otimes \hat{m}^{\op}"']&B\ar[d, "\hat{m}^{\op}"]\\
		\Map(B,B)\otimes\Map(B,B)\ar[d, "C"']&\Map(B,B)\\
		\Map(B,B)\otimes\Map(B,B).\ar[ru, "\circ"']
		\end{tikzcd}\]
		Since the enriched composition is associative (\cite{MR2177301} 1.6), the assertion follows by comparing the two diagrams
		\[\begin{tikzcd}
		B\otimes B\otimes \Map(B,M)\ar[r, "\hat{m}^{\op}"]\ar[d, "\hat{m}^{\op}\otimes \hat{m}^{\op}"']
		&B\otimes\Map(B,B)\otimes\Map(B,M)\ar[d, "\id_B\otimes C"]\ar[ld, "\hat{m}^{\op}"]\\
		\Map(B,B)\otimes\Map(B,B)\otimes\Map(B,M)\ar[dd, "\id_{\Map(B,B)}\otimes C"']&B\otimes \Map(B,M)\otimes \Map(B,B)\ar[d, "\circ"]\\
		&B\otimes \Map(B,M)\ar[d, "\hat{m}^{\op}"]\\
		\Map(B,B)\otimes\Map(B,M)\otimes\Map(B,B)\ar[r, "\id\otimes\circ"]\ar[d, "C_{\Map(B,B),\Map(B,M)\otimes \Map(B,B)}"]
		&\Map(B,B)\otimes\Map(B,M)\ar[d, "C"]\\
		\Map(B,M)\otimes\Map(B,B)\otimes\Map(B,B)\ar[r, "\circ\otimes \id"]&\Map(B,M)\otimes\Map(B,B)\ar[d, "\circ"]\\
		&\Map(B,M)
		\end{tikzcd}\]
		\[\begin{tikzcd}
		B\otimes B\otimes \Map(B,M)\ar[r, "m"]\ar[d, "\hat{m}^{\op}\otimes \hat{m}^{\op}"']
		&B\otimes\Map(B,M)\ar[d, "\hat{m}^{\op}"]\\
		\Map(B,B)\otimes\Map(B,B)\otimes\Map(B,M)\ar[d, "C\otimes \id_{\Map(B,M)}"']&\Map(B,B)\otimes\Map(B,M)\ar[d, "C"]\\
		\Map(B,B)\otimes\Map(B,B)\otimes\Map(B,M)\ar[ru, "\circ\otimes \id_{\Map(B,M)}"']\ar[d, "C_{\Map(B,B)\otimes\Map(B,B),\Map(B,M)}"']
		&\Map(B,M)\otimes\Map(B,B)\ar[d, "\circ"]\\
		\Map(B,M)\otimes\Map(B,B)\otimes\Map(B,B)\ar[ru, "\id_{\Map(B,M)}\otimes \circ"']&\Map(B,M).
		\end{tikzcd}\]
		For the coincidence of the left vertical composite arrows, see \cite{MR1250465} Proposition 2.7 and also \cite{JS} Proposition 1 B.5. Observe next that the map $\hat{m}^{\op}\circ j_B:I\to \Map(B,B)$ corresponds to $l:I\otimes B\cong B$ (use \cite{MR1250465} Proposition 2.1). Then the unitality follows from the following diagram by passing to the left adjoints:
		\[\begin{tikzcd}
		I\otimes\Map(B,M)\otimes B\ar[d, "C\otimes \id_B"']\ar[rrd, "l_{\Map(B,M)}"]&\\
		\Map(B,M)\otimes I\otimes B\ar[d, "l_B"']\ar[rr, "r_{\Map(B,M)}"]&&\Map(B,M)\otimes B\ar[lld, equal]\\
		\Map(B,M)\otimes B\ar[r, "\ev"']&M.
		\end{tikzcd}\]
		
		We next define two morphisms $\Map(B,M)\rightrightarrows\Map(A\otimes B,M)$ as follows:
		\begin{itemize}
			\item Take the pullback along $m_B\circ (f\otimes \id_B):A\otimes B\to B\otimes B\to B$.
			\item Define $\Map(B,M)\to\Map(A\otimes B,A\otimes M)$ by
			\[\Map(B,M)\otimes A\otimes B\overset{C}{\cong}A\otimes\Map(B,M)\otimes B\overset{\ev}{\to}A\otimes M.\]
			Then compose it with $\pi_M$.
		\end{itemize}
		Let $\Map_A(B,M)$ be the equalizer of these two maps, and $i:\Map_A(B,M)\to\Map(B,M)$ denote the canonical morphism. Then the composite two arrows
		\[B\otimes\Map_A(B,M)\overset{i}{\to} B\otimes\Map(B,M)\overset{\pi}{\to}\Map(B,M)\rightrightarrows \Map(A\otimes B,M)\]
		coincide. In particular, $\Map_A(B,M)$ is a $B$-submodule of $\Map(B,M)$. In fact, passing to left adjoints, we can rewrite the composite arrows $B\otimes\Map_A(B,M)\rightrightarrows \Map(A\otimes B,M)$ as
		\[\begin{split}
		B\otimes\Map_A(B,M)\otimes A\otimes B&\overset{i}{\to} B\otimes\Map(B,M)\otimes A\otimes B\\
		&\overset{\hat{m}^{\op}}{\to}\Map(B,B)\otimes\Map(B,M)\otimes A\otimes B\\
		&\overset{C}{\cong}\Map(B,M)\otimes\Map(B,B)\otimes A\otimes B\\
		&\overset{\circ}{\to}\Map(B,M)\otimes A\otimes B\\
		&\overset{f}{\to}\Map(B,M)\otimes B\otimes B\\
		&\overset{m}{\to}\Map(B,M)\otimes B\\
		&\overset{\ev}{\to}M
		\end{split}\]
		\[\begin{split}
		B\otimes\Map_A(B,M)\otimes A\otimes B&\overset{i}{\to} B\otimes\Map(B,M)\otimes A\otimes B\\
		&\overset{\hat{m}^{\op}}{\to}\Map(B,B)\otimes\Map(B,M)\otimes A\otimes B\\
		&\overset{C}{\cong}\Map(B,M)\otimes\Map(B,B)\otimes A\otimes B\\
		&\overset{\circ}{\to}\Map(B,M)\otimes A\otimes B\\
		&\overset{C}{\cong}A\otimes \Map(B,M)\otimes B\\
		&\overset{\ev}{\to}A\otimes M\\
		&\overset{\pi_M}{\to} M.
		\end{split}\]
		Since the arrows $\Map_A(B,M)\overset{i}{\to}\Map(B,M)\rightrightarrows\Map(A\otimes B,M)$ coincide, the diagram
		\[\begin{tikzcd}
		\Map_A(B,M)\otimes A\otimes B\ar[r, "i"]\ar[d, "i"']&\Map(B,M)\otimes A\otimes B\ar[d, "C\otimes\id_B"]\\
		\Map(B,M)\otimes A\otimes B\ar[d, "f"']&A\otimes\Map(B,M)\otimes B\ar[d, "\ev"]\\
		\Map(B,M)\otimes B\otimes B\ar[d, "m_B"']&A\otimes M\ar[d, "\pi_M"]\\
		\Map(B,M)\otimes B\ar[r, "\ev"]&M
		\end{tikzcd}\]
		commutes. The assertion now follows by comparing diagrams
		\[\begin{tikzcd}
		B\otimes\Map_A(B,M)\otimes A\otimes B\ar[r, "C\otimes \id_{A\otimes B}"]
		\ar[d, "\hat{m}^{\op}\otimes i\otimes\id_{A\otimes B}"']
		&\Map_A(B,M)\otimes B\otimes A\otimes B\ar[dd, "i\otimes\hat{m}^{\op}\otimes f\otimes\id_B"]\\
		\Map(B,B)\otimes \Map(B,M)\otimes A\otimes B\ar[d, "C"']&\\
		\Map(B,M)\otimes\Map(B,B)\otimes A\otimes B\ar[d, "\circ"']\ar[r, "f"']&\Map(B,M)\otimes\Map(B,B)\otimes B\otimes B\ar[d, "m"]\\
		\Map(B,M)\otimes A\otimes B\ar[d, "m_B\circ (f\otimes \id_B)"']&\Map(B,M)\otimes\Map(B,B)\otimes B\ar[ld, "\circ"]\ar[d, "\ev"]\\
		\Map(B,M)\otimes B\ar[d, "\ev"']&\Map(B,M)\otimes B\ar[ld, "\ev"]\\
		M
		\end{tikzcd}\]
		\[\begin{tikzcd}
		\Map_A(B,M)\otimes B\otimes A\otimes B\ar[r, "C_{B,A\otimes B}"]\ar[d, "i\otimes\hat{m}^{\op}\otimes f"']
		&\Map_A(B,M)\otimes A\otimes B\otimes B\ar[d, "f"]\\
		\Map(B,M)\otimes\Map(B,B)\otimes B\otimes B\ar[d, "m"']&\Map_A(B,M)\otimes B\otimes B\otimes B
		\ar[d, "i\otimes m\circ(m\otimes \id_B)"]\\
		\Map(B,M)\otimes\Map(B,B)\otimes B\ar[r, "\ev"']&\Map(B,M)\otimes B
		\end{tikzcd}\]
		\[\begin{tikzcd}
		B\otimes\Map_A(B,M)\otimes A\otimes B\ar[d, "\hat{m}^{\op}\otimes i\otimes \id_A"']\ar[r, "C_{B,\Map_A(B,M)\otimes A}"]
		&\Map_A(B,M)\otimes A\otimes B\otimes B\ar[d, "i\otimes\id_A\otimes\hat{m}^{\op}"]\\
		\Map(B,B)\otimes \Map(B,M)\otimes A\otimes B\ar[d, "C"']\ar[r, "C"]&\Map(B,M)\otimes A\otimes\Map(B,B)\otimes B\ar[d, "C"]\\
		\Map(B,M)\otimes\Map(B,B)\otimes A\otimes B\ar[d, "\circ"']\ar[r, "C"]
		&A\otimes\Map(B,M)\otimes\Map(B,B)\otimes B\ar[d, "\ev"]\ar[ldd, "\circ"']\\
		\Map(B,M)\otimes A\otimes B\ar[d, "C"']&A\otimes\Map(B,M)\otimes B\ar[d, "\ev"]\\
		A\otimes\Map(B,M)\otimes B\ar[r, "\ev"']&A\otimes M\ar[d, "\pi_M"]\\
		&M
		\end{tikzcd}\]
		\[\begin{tikzcd}
		\Map_A(B,M)\otimes A\otimes B\otimes B\ar[d, "i\otimes\id_A\otimes\hat{m}^{\op}"']\ar[r, "i"]
		&\Map(B,M)\otimes A\otimes B\otimes B\ar[dd, "C_{\Map(B,M),A}"]\ar[ld, "\hat{m}^{\op}"]\\
		\Map(B,M)\otimes A\otimes\Map(B,B)\otimes B\ar[d, "C_{\Map(B,M),A}"']&\\
		A\otimes\Map(B,M)\otimes\Map(B,B)\otimes B\ar[d, "\ev"']&A\otimes\Map(B,M)\otimes B\otimes B\ar[l, "\hat{m}^{\op}"]\ar[d, "C_{B,B}"]\\
		A\otimes\Map(B,M)\otimes B\ar[d, "\ev"']&A\otimes\Map(B,M)\otimes B\otimes B\ar[l, "m"']\\
		A\otimes M.
		\end{tikzcd}\]
		To see that this defines a right adjoint functor to the forgetful functor $B\cmod \to A\cmod$, let $M$ be a $B$-module, and $N$ be an $A$-module. If we are given an $A$-module homomorphism $p:M\to N$, define a $B$-module homomorphism $M\to\Map (B,N)$ adjunctionally by
		\[M\otimes B\overset{C_{M,B}}{\cong} B\otimes M\overset{\pi}{\to}M\overset{p}{\to} N.\]
		This factors through $\Map_A(B,M)$ by a similar argument. For a $B$-module homomorphism $M\to \Map_A(B,N)$, define an $A$-module homomorphism by
		\[M\to\Map_A(B,N)\to\Map(B,N)\to N,\]
		where the last arrow is defined by restriction to $I$: $\Map(B,N)\to\Map(I,N)\cong N$.
	\end{proof}
	\begin{lem}\label{2.2.7}
		Let $\cV$ be a symmetric monoidal category with coequalizers, and $f:(A,m_A)\to (B,m_B)$ be a morphism of monoids of $\cV$. Suppose that for every object $X\in\cV$, $X\otimes -$ respects coequalizers. Then the forgetful functor
		\[B\cmod\to A\cmod\]
		admits a left adjoint functor.
	\end{lem}
	\begin{proof}
		The proof essentially goes in the dual way to Lemma \ref{2.2.6}. Let $(M,\pi_M)$ be an $A$-module. It is clear that $B\otimes M$ is a left $B$-module for the multiplication from the left side. Define $B\otimes_A M$ by the coequalizer sequence
		\[B\otimes A\otimes M\xrightrightarrows[\id_B\otimes \pi_M]{(m_B\circ (\id_B\otimes f))\otimes \id_M} B\otimes M\overset{q}{\to} B\otimes_A M.\]
		Then we have a commutative diagram
		\[\begin{tikzcd}
		&B\otimes B\otimes A\otimes M\ar[rd, "m_B"']\ar[r, "\pi_M"]\ar[d, "f"']
		&B\otimes B\otimes M\ar[r, "m_B\otimes\id_M"]&B\otimes M\ar[ddd, "q"]\\
		&B\otimes B\otimes B\otimes M\ar[d, "m_B\otimes \id_B"']\ar[ld, "\id_B\otimes m_B"']
		&B\otimes A\otimes M\ar[ru, "\id_B\otimes\pi_M"']\ar[ldd, "(m_B\circ (\id_B\otimes f))\otimes\id_M"]&&\\
		B\otimes B\otimes M\ar[rd, "m_B"']&B\otimes B\otimes M\ar[d, "m_B"']&&&\\
		&B\otimes M\ar[rr, "q"']&&B\otimes_A M.&
		\end{tikzcd}\]
		Since $B\otimes-$ respects coequalizer sequences, this implies that the action of $B$ on $B\otimes M$ descends to $B\otimes_A M$. Moreover, the resulting map $B\otimes (B\otimes_A M)\to B\otimes_A M$ exhibits the structure of a left $B$-module on $B\otimes_A M$ since $(B\otimes B)\otimes-$ respects coequalizers sequences. The adjunction is defined in the usual manner: Let $M$ be an $A$-module, and $N$ be a $B$-module. Suppose that we are given an $A$-module homomorphism $p:M\to N$. Then it extends to a $B$-module homomorphism $\pi_N\circ (\id_B\otimes p):B\otimes M\to N$. A similar argument shows that this extension of $p$ descends to a $B$-module homomorphism $B\otimes_A M\to N$:
		\[\begin{tikzcd}
		B\otimes A\otimes M\ar[rr, "\pi_M"]\ar[rd, "p"']\ar[d, "f"']&&B\otimes M\ar[dd, "p"]\\
		B\otimes B\otimes M\ar[rd, "p"']\ar[d, "m_B"']&B\otimes A\otimes N\ar[rd, "\pi_N"]\ar[d, "f"']\\
		B\otimes M\ar[d, "p"']&B\otimes B\otimes N\ar[r, "\pi_N"']\ar[ld, "m_B"']&B\otimes N\ar[d, "\pi_N"]\\
		B\otimes N\ar[rr, "\pi_N"']&&N
		\end{tikzcd}\]
		This determines the desired adjunction.
	\end{proof}
	\begin{thm}\label{2.2.8}
		\begin{enumerate}
			\renewcommand{\labelenumi}{(\arabic{enumi})}
			\item For a map $(\cA,K)\to(\cB,L)$ of pairs, $\cF_{\cB,L}^{\cA,K}$ admits a right adjoint functor $I^{\cB,L}_{\cA,K}$.
			\item For a map $f=(f_a,\id_K):(\cA,K)\to (\cB,K)$ of pairs, the forgetful functor $\cF_{\cB,K}^{\cA,K}$ admits a left adjoint functor $P^{\cB,K}_{\cA,K}$.
		\end{enumerate}
	\end{thm}
	\begin{proof}
		Recall that a given map $(\cA,K)\to (\cB,L)$ of (weak) pairs admits a factorization
		\[(\cA,K)\to (\cB,K)\to(\cB,L)\]
		(see above Proposition \ref{2.2.5}), which implies an equality
		\[\cF^{\cA,K}_{\cB,L,w}=\cF_{\cB,K,w}^{\cA,K}\circ\cF^{\cB,K}_{\cB,L,w}.\]
		Notice that $\cF^{k,K}_{k,L,w}$ is a (strict) symmetric monoidal functor. Moreover, it admits a right adjoint functor (\cite{MR2015057} I.3.3 and I.3.4). We now apply Lemma \ref{2.2.5} to 
		\[\cV=L\cmod;\]
		\[\cV'=K\cmod;\]
		\[F=\cF^{k,K}_{k,L,w};\]
		\[A=(\cB,K)\]
		to get a right adjoint functor $I^{\cB,L}_{\cB,K,w}$ of $\cF^{\cB,K}_{\cB,L,w}$. Similarly, apply
		\[\cV=K\cmod\]
		\[A=(\cA,K)\]
		\[B=(\cB,K)\]
		to Lemma \ref{2.2.6} and Lemma \ref{2.2.7} to get right and left adjoint functors $I^{\cB,K}_{\cA,K,w}$ and $P^{\cB,K}_{\cA,K,w}$. Compose the right adjoint functors to get a right adjoint functor $I^{\cB,L}_{\cA,K,w}$ of $\cF^{\cA,K}_{\cB,L,w}$. From Lemma \ref{2.2.2}, 
		\[I^{\cB,L}_{\cA,K}=(-)^\fk\circ I^{\cB,L}_{\cA,K,w}\circ\cJ_{\cA,K}\]
		is right adjoint to $\cF^{\cA,K}_{\cB,L}$. This shows (1).
		
		To see (2), it will suffice to show that $P^{\cB,K}_{\cA,K,w}$ sends $(\cA,K)$-modules to $(\cB,K)$-modules. Let $(M,\pi_M,\nu_M)$ be an $(\cA,K)$-module and set $X:=P^{\cB,K}_{\cA,K,w}M$. Let us denote the corresponding action maps on $X$ by $\pi_X$ and $\nu_X$. Let $\psi_\cA:\fk\to\cA$ and $\psi_\cB:\fk\to\cB$ be the structure maps. Then for $\xi\in\fk$ and $b\otimes m\in X$, we obtain
		\[\begin{split}d\nu_X(\xi)(b\otimes m)&=\left[\psi_\cB(\xi),b\right]\otimes m+b\otimes d\nu_M(\xi)m\\
		&=\psi_\cB(\xi)b\otimes m-b\otimes\pi_M(\psi_\cA(\xi))m+b\otimes d\nu_M(\xi)m\\
		&=\psi_\cB(\xi)b\otimes m\\
		&=\pi_X(\psi_\cB(\xi))(b\otimes m).
		\end{split}\]
		This completes the proof.
	\end{proof}
	\begin{rem}\label{2.2.9}
		The functor $\cF_{\cB,L}^{\cA,K}$ does not always admit a left adjoint functor since it does not preserve infinite limits in general.
	\end{rem}
	\begin{rem}[production-in-stages]\label{2.2.10}
		Suppose that we are given a sequence
		\[(\cA,K)\to (\cB,L) \to (\cC,M)\]
		of pairs. Then we have an obvious equality
		\[\cF_{\cB,L}^{\cA,K}\circ\cF^{\cB,L}_{\cC,M}=\cF^{\cA,K}_{\cC,M}.\]
		Pass to the right adjoints to get an isomorphism $I_{\cB,L}^{\cC,M}\circ I^{\cB,L}_{\cA,K}\cong I_{\cA,K}^{\cC,M}$.
	\end{rem}
	Finally, we put a linear structure on the category of weak $(\cA,K)$-modules. Notice that each Hom set of the category of weak $(\cA,K)$-modules is a $k$-submodule of the $k$-module of morphisms as $k$-modules. The resulting linear structure makes $(\cA,K)\cmod_w$ into a linear category. Remark that the isomorphism of Lemma \ref{2.1.7} preserves the linear structures. The next result is a consequence of Lemma \ref{2.2.2}:
	\begin{cor}\label{2.2.11}
		\begin{enumerate}
			\renewcommand{\labelenumi}{(\arabic{enumi})}
			\item For a weak pair $(\cA,K)$, the category $(\cA,K)\cmod_w$ is a (locally small) bicomplete abelian category, whose colimits and finite limits are computed in the category of dg $k$-modules.
			\item For a pair $(\cA,K)$, the category $(\cA,K)\cmod$ is a (locally small) bicomplete abelian category, whose colimits and finite limits are computed in the category of dg $k$-modules.
		\end{enumerate}
	\end{cor}
	\begin{proof}
		Let $(\cA,K)$ be a weak pair. According to \cite{MR1313497} Proposition 4.3.1, Proposition 4.3.2, $(\cA,K)\cmod_w$ has small limits and colimits which are computed in $K\cmod$ since $K\cmod$ is a bicomplete closed monoidal category. Moreover, colimits and finite limits are computed in $k\cmod$ since $K$ is flat over the base ring $k$. Notice that the zero dg $k$-module admits the trivial weak $(\cA,K)$-module structure, which is a zero object of $(\cA,K)\cmod_w$. It is also clear that $(\cA,K)\cmod_w$ is additive. Finally, suppose that we are given a map of weak $(\cA,K)$-modules $f:M\to N$. Let $\Ker f$ (resp.\ $\Image f$) denote the kernel of $f$ (resp.\ the image of $f$). Then we have a natural map
		\[M/\Ker f\to\Image f\]
		which is an isomorphism as a map of $k\cmod$, so it is also an isomorphism as a map of weak $(\cA,K)$-modules. Hence (1) is verified.
		
		We next prove (2). Let $(\cA,K)$ be a pair. Lemma \ref{2.2.2} implies that the category $(\cA,K)\cmod$ is bicomplete and that its limits and colimits are computed in $(\cA,K)\cmod_w$. In particular, $(\cA,K)\cmod$ is an abelian subcategory of $(\cA,K)\cmod$. This completes the proof.
	\end{proof}
	\begin{rem}\label{2.2.12}
		We can directly prove without Lemma \ref{2.2.2} that the category $(\cA,K)\cmod$ is stable under formation of (possibly empty or infinite) coproducts, cokernels, and kernels in $(\cA,K)\cmod_w$ by checking the compatibility condition of the actions of $\fk$.
	\end{rem}
	\section{Differential graded analogue of the dual Zuckerman functor}
	In this section, we remark how to define differential graded enrichment of the dual Zuckerman functor and the pseudoforgetful functor over the field $\bC$ of complex numbers. In this section, the ground field is $\bC$, and all affine group schemes are assumed to be reductive. We work with reductive groups since their representations are completely reducible. Thanks to this fact, the pseudoforgetful functor is exact (\cite{MR1330919} Lemma 2.28 and Proposition 2.33). We work over $\bC$ to use \cite{MR1330919} Lemma 1.46. The author expects that \cite{MR1330919} Lemma 1.46 is valid in the setting of \cite{1604.04253}.
	\subsection{Construction}
	For a reductive group $K$, $R(K)$ will denote the Hecke algebra of a maximal compact subgroup of $K$ (\cite{MR1330919} I.2). Recall that $R(K)$ is a right $\cO(K)$-module by the standard multiplication (\cite{MR1330919} Lemma 1.46).
	
	As in the proof of Theorem \ref{2.2.8} (1), we start with the weak setting. Let $(\cA,L,\phi)$ be a weak pair, $K\to L$ be a homomorphism of reductive groups, and $V$ be a weak $(\cA,K)$-module. As a dg $L$-module, define $P^{\cA,L}_{\cA,K,w}(M)=\Pi_w(M)$ as
	\[\Pi_w(M)=R(L)\otimes_{R(K)} M.\]
	To describe the action of $\cA$ on $\Pi_w(M)$, let $a\in\cA$, and $T\otimes m\in \Pi_w(M)$. Choose elements $a_i\in\cA$ and $f_i\in\cO(K)$ such that
	\[\phi(l)^{-1}a=\sum_i f_i(l)a_i\]
	for every $l\in L$. If we write $\rho$ for the coaction of $\cO(L)$ on $\cA$, these elements are characterized by
	\[\rho(a)=\sum_i a_i\otimes S(f_i),\]
	where $S$ denotes the antipode of $\cO(L)$. Then the action $\pi_{\Pi_w(M)}(a)$ on $\Pi_w(M)$ is written as
	\[\pi_{\Pi_w(M)}(a)(T\otimes m)=\sum Tf_i\otimes \pi(a_i)m.\]
	This is clearly independent of the expression of $\rho(a)$, and determines a differential graded action (\cite{MR1330919} I.5, \cite{MR2244116} 5.4.1). The dual Zuckerman functor for pairs is defined by composing $\Pi_w$ with the localization $(-)_\fl$.
	
	Similarly, let $N$ be a weak $(\cA,L)$-module. Define a dg vector subspace $\Map_{R(L)}(R(L),N)\subset\Map(R(L),N)$ by
	\[\Map_{R(L)}(R(L),N)^i=\{\mathrm{graded\ left\ }R(L)
	\textrm{-}\mathrm{module\ homomorphisms\ }
	R(L)\to N\left[i\right]\}.\]
	Then $\Map_{R(L)}(R(L),N)$ is a left $R(K)$-module for the right action on the domain. As a $K$-module, $(\cF^\vee)^{\cA,K}_{\cA,L,w}(N)$ is defined as its $K$-finite part. To define the action of $\cA$, take homogeneous elements $\varphi\in\Map_{R(L)}(R(L),M)_K$, $a\in\cA$ and $S\in R(L)$. Choose $a_i\in\cA$ and $f_i\in\cO(L)$ so that
	\[\phi(l)a=\sum f_i(l)a_i\]
	for every $l\in L$. Then the action $\pi(a)=\pi_{(\cF^\vee)^{\cA,K}_{\cA,L,w}(N)}(a)$ is defined as
	\[(\pi(a)\varphi)(T)=\sum \pi_M(a_i)\varphi(Tf_i).\]
	This determines the structure of a weak $(\cA,K)$-module. Since the action of $\fk$ induced from $\cA$ factors through the degree 0 part $\cA^0$, this functor sends $(\cA,L)$-modules to $(\cA,K)$-modules if $(\cA,L)$ is a pair from the classical setting. We define $(\cF^{\vee})_{\cA ,L}^{\cA ,K}: (\cA,L)\cmod \to (\cA,K)\cmod$ by the restriction of $(\cF^{\vee})^{\cA,K}_{\cA,L,w}$ if $(\cA,L)$ is a pair.
	\subsection{Induction-in-stages}
	To establish the principle of induction-in-stages, we extend our definition. If we are given a map $f:(\cA,K)\to(\cB,L)$ of pairs, factorize it as
	\[\begin{tikzcd}
	(\cA,K)\ar[rr, "f"]\ar[rd, "(f_a\mathrm{,} \id_K)"']&&(\cB,L)\\
	&(\cB,K).\ar[ru, "(\id_{\cB}\mathrm{,}f_k)"']
	\end{tikzcd}\]
	According to this factorization, set
	\[P^{\cB,L}_{\cA,K}=P^{\cB,L}_{\cB,K}\circ P^{\cB,K}_{\cA,K};\]
	\[(\cF^\vee)_{\cB,L}^{\cA,K}=\cF_{\cB,K}^{\cA,K}\circ(\cF^\vee)_{\cB,L}^{\cB,K}.\]
	\begin{prop}[induction-in-stages]
		Suppose that we are given a sequence $(\cA,K)\to(\cB,L)\to(\cC,M)$ of maps of pairs. Then the following isomorphisms exist:
		\begin{enumerate}
			\renewcommand{\labelenumi}{(\arabic{enumi})}
			\item $(\cF^\vee)_{\cB,L}^{\cA,K}\circ(\cF^\vee)^{\cB,L}_{\cC,M}\cong (\cF^\vee)^{\cA,K}_{\cC,M}$;
			\item $P_{\cB,L}^{\cC,M}\circ P^{\cB,L}_{\cA,K}\cong P_{\cA,K}^{\cC,M}$.
		\end{enumerate}
	\end{prop}
	\begin{proof}
		It suffices to prove the weak version of (1). We construct the isomorphism in the following two stages:
		\[\begin{split}
		(\cF^\vee)_{\cB,L,w}^{\cA,K}\circ(\cF^\vee)^{\cB,L}_{\cC,M,w}
		&=\cF_{\cB,K,w}^{\cA,K}\circ(\cF^\vee)_{\cB,L,w}^{\cB,K}\circ\cF_{\cC,L,w}^{\cB,L}\circ(\cF^\vee)^{\cC,L}_{\cC,M,w}\\
		&\stackrel{(i)}{\cong}\cF_{\cB,K,w}^{\cA,K}\circ\cF_{\cC,K,w}^{\cB,K}\circ(\cF^\vee)_{\cC,L,w}^{\cC,K}\circ(\cF^\vee)_{\cC,M,w}^{\cC,L}\\
		&\stackrel{(ii)}{\cong}\cF_{\cB,K,w}^{\cA,K}\circ\cF_{\cC,K,w}^{\cB,K}\circ
		(\cF^\vee)_{\cC,M,w}^{\cC,K}
		\\
		&=\cF_{\cC,K,w}^{\cA,K}\circ
		(\cF^\vee)_{\cC,M,w}^{\cC,K}\\
		&=(\cF^\vee)^{\cA,K}_{\cC,M,w}.
		\end{split}\]
		Part (i) follows by unwinding the definitions. To see (ii), we pass to the left adjoints again. Since the natural isomorphism
		\[R(M)\otimes_{R(L)}(R(L)\otimes_{R(K)}-)\cong R(M)\otimes_{R(K)}-\]
		respects the differential and the actions of $\cC$ and $R(M)$, the assertion follows. This completes the proof.
	\end{proof}
	

\begin{thebibliography}{10}
		
		\bibitem{MR2724388}
		M.~Aguiar and S.~Mahajan.
		\newblock {\em Monoidal functors, species and {H}opf algebras}, volume~29 of
		{\em CRM Monograph Series}.
		\newblock American Mathematical Society, Providence, RI, 2010.
		\newblock With forewords by Kenneth Brown and Stephen Chase and Andr\'{e}
		Joyal.
		
		\bibitem{MR1237825}
		A.~Beilinson and J.~Bernstein.
		\newblock A proof of {J}antzen conjectures.
		\newblock In {\em I. {M}. {G}el'fand {S}eminar}, volume~16 of {\em Adv. Soviet
			Math.}, pages 1--50. Amer. Math. Soc., Providence, RI, 1993.
		
		\bibitem{10.1093/imrn/rny147}
		J.~Bernstein, N.~Higson, and E.~Subag.
		\newblock {Algebraic Families of Harish-Chandra Pairs}.
		\newblock {\em Int. Math. Res. Not. IMRN}, 07 2018.
		
		\bibitem{10.1093/imrn/rny146}
		J.~Bernstein, N.~Higson, and E.~Subag.
		\newblock {Contractions of Representations and Algebraic Families of
			Harish-Chandra Modules}.
		\newblock {\em Int. Math. Res. Not. IMRN}, 07 2018.
		
		\bibitem{MR1317229}
		J.~Bernstein and V.~Lunts.
		\newblock Localization for derived categories of {$({\mathfrak g},K)$}-modules.
		\newblock {\em J. Amer. Math. Soc.}, 8(4):819--856, 1995.
		
		\bibitem{MR1313497}
		F.~Borceux.
		\newblock {\em Handbook of categorical algebra. 2}, volume~51 of {\em
			Encyclopedia of Mathematics and its Applications}.
		\newblock Cambridge University Press, Cambridge, 1994.
		\newblock Categories and structures.
		
		\bibitem{MR563524}
		M.~Demazure and P.~Gabriel.
		\newblock {\em Introduction to algebraic geometry and algebraic groups},
		volume~39 of {\em North-Holland Mathematics Studies}.
		\newblock North-Holland Publishing Co., Amsterdam-New York, 1980.
		\newblock Translated from the French by J. Bell.
		
		\bibitem{MR2138086}
		L.~El~Kaoutit, J.~G\'{o}mez-Torrecillas, and F.~J. Lobillo.
		\newblock Semisimple corings.
		\newblock {\em Algebra Colloq.}, 11(4):427--442, 2004.
		
		\bibitem{MR925070}
		V.~A. Ginzburg.
		\newblock Equivariant cohomology and {K}\"{a}hler geometry.
		\newblock {\em Funktsional. Anal. i Prilozhen.}, 21(4):19--34, 96, 1987.
		
		\bibitem{1407.0574}
		G.~Harder.
		\newblock {H}arish-{C}handra {M}odules over {$\mathbb Z$}, 2014.
		\newblock arXiv:1405.6513.
		
		\bibitem{Harder2011349}
		G.~Harder and A.~Raghuram.
		\newblock Eisenstein cohomology for {GL}$_{N}$ and the special values of
		{R}ankin-{S}elberg {$L$}-functions.
		\newblock 2019.
		\newblock To appear as a volume in the Annals of Math. Studies.
		
		\bibitem{MR3053412}
		M.~Harris.
		\newblock Beilinson-{B}ernstein localization over {$\mathbb Q$} and periods of
		automorphic forms.
		\newblock {\em Int. Math. Res. Not. IMRN}, (9):2000--2053, 2013.
		
		\bibitem{10.1093/imrn/rny043}
		M.~Harris.
		\newblock Beilinson-{B}ernstein localization over {$\mathbb Q$} and periods of
		automorphic forms: Erratum.
		\newblock {\em Int. Math. Res. Not. IMRN}, 03 2018.
		
		\bibitem{1606.04320}
		T.~Hayashi.
		\newblock Dg analogues of the {Z}uckerman functors and the dual {Z}uckerman
		functors {II}, 2016.
		\newblock arXiv:1606.04320.
		
		\bibitem{1712.07336}
		T.~Hayashi.
		\newblock Integral models of {H}arish-{C}handra modules of the finite covering
		groups of {PU}(1,1), 2017.
		\newblock arXiv:1712.07336.
		
		\bibitem{MR3853058}
		T.~Hayashi.
		\newblock Flat base change formulas for {$(\mathfrak{g},K)$}-modules over
		{N}oetherian rings.
		\newblock {\em J. Algebra}, 514:40--75, 2018.
		
		\bibitem{MR2066503}
		M.~Hovey.
		\newblock Homotopy theory of comodules over a {H}opf algebroid.
		\newblock In {\em Homotopy theory: relations with algebraic geometry, group
			cohomology, and algebraic {$K$}-theory}, volume 346 of {\em Contemp. Math.},
		pages 261--304. Amer. Math. Soc., Providence, RI, 2004.
		
		\bibitem{MR2244116}
		J-S. Huang and P.~Pand\v{z}i\'{c}.
		\newblock {\em Dirac operators in representation theory}.
		\newblock Mathematics: Theory \& Applications. Birkh\"{a}user Boston, Inc.,
		Boston, MA, 2006.
		
		\bibitem{MR2015057}
		J.~C. Jantzen.
		\newblock {\em Representations of algebraic groups}, volume 107 of {\em
			Mathematical Surveys and Monographs}.
		\newblock American Mathematical Society, Providence, RI, second edition, 2003.
		
		\bibitem{1604.04253}
		F.~Januszewski.
		\newblock On {P}eriod {R}elations for {A}utomorphic {$L$}-functions {II}, 2016.
		\newblock arXiv:1604.04253.
		
		\bibitem{1808.10709}
		F.~Januszewski.
		\newblock Families of {$\mathcal D$}-modules, 2018.
		\newblock arXiv:1808.10709.
		
		\bibitem{MR3770183}
		F.~Januszewski.
		\newblock Rational structures on automorphic representations.
		\newblock {\em Math. Ann.}, 370(3-4):1805--1881, 2018.
		
		\bibitem{MR3937337}
		F.~Januszewski.
		\newblock On period relations for automorphic {$L$}-functions {I}.
		\newblock {\em Trans. Amer. Math. Soc.}, 371(9):6547--6580, 2019.
		
		\bibitem{JS}
		A.~Joyal and R.~Street.
		\newblock {\em Braided monoidal categories}.
		\newblock Macquarie Math. Reports No. 860081, 1986.
		
		\bibitem{MR1250465}
		A.~Joyal and R.~Street.
		\newblock Braided tensor categories.
		\newblock {\em Adv. Math.}, 102(1):20--78, 1993.
		
		\bibitem{MR2177301}
		G.~M. Kelly.
		\newblock Basic concepts of enriched category theory.
		\newblock {\em Repr. Theory Appl. Categ.}, (10):vi+137, 2005.
		\newblock Reprint of the 1982 original [Cambridge Univ. Press, Cambridge;
		MR0651714].
		
		\bibitem{MR2945222}
		S.~N. Kitchen.
		\newblock Cohomology of standard modules on partial flag varieties.
		\newblock {\em Represent. Theory}, 16:317--344, 2012.
		
		\bibitem{MR1330919}
		A.~W. Knapp and D.~A. Vogan, Jr.
		\newblock {\em Cohomological induction and unitary representations}, volume~45
		of {\em Princeton Mathematical Series}.
		\newblock Princeton University Press, Princeton, NJ, 1995.
		
		\bibitem{luriehigheralgebra}
		J.~Lurie.
		\newblock Higher {A}lgebra.
		\newblock \url{http://www.math.harvard.edu/ lurie/}.
		
		\bibitem{MR2522659}
		J.~Lurie.
		\newblock {\em Higher topos theory}, volume 170 of {\em Annals of Mathematics
			Studies}.
		\newblock Princeton University Press, Princeton, NJ, 2009.
		
		\bibitem{MR1712872}
		S.~Mac~Lane.
		\newblock {\em Categories for the working mathematician}, volume~5 of {\em
			Graduate Texts in Mathematics}.
		\newblock Springer-Verlag, New York, second edition, 1998.
		
		\bibitem{MR1486143}
		D.~Mili\v{c}i\'{c} and P.~Pand\v{z}i\'{c}.
		\newblock Equivariant derived categories, {Z}uckerman functors and
		localization.
		\newblock In {\em Geometry and representation theory of real and {$p$}-adic
			groups ({C}\'{o}rdoba, 1995)}, volume 158 of {\em Progr. Math.}, pages
		209--242. Birkh\"{a}user Boston, Boston, MA, 1998.
		
		\bibitem{MR2692966}
		P.~Pand\v{z}i\'{c}.
		\newblock {\em Equivariant analogues of {Z}uckerman functors}.
		\newblock ProQuest LLC, Ann Arbor, MI, 1995.
		\newblock Thesis (Ph.D.)--The University of Utah.
		
		\bibitem{MR2171292}
		P.~Pand\v{z}i\'{c}.
		\newblock A simple proof of {B}ernstein-{L}unts equivalence.
		\newblock {\em Manuscripta Math.}, 118(1):71--84, 2005.
		
		\bibitem{MR2276617}
		P.~Pand\v{z}i\'{c}.
		\newblock Zuckerman functors between equivariant derived categories.
		\newblock {\em Trans. Amer. Math. Soc.}, 359(5):2191--2220, 2007.
		
		\bibitem{MR1680015}
		J.~C. Santos.
		\newblock Foncteurs de {Z}uckerman pour les superalg\`ebres de {L}ie.
		\newblock {\em J. Lie Theory}, 9(1):69--112, 1999.
		
		\bibitem{MR3827131}
		E.~M. Subag.
		\newblock Symmetries of the hydrogen atom and algebraic families.
		\newblock {\em J. Math. Phys.}, 59(7):071702, 20, 2018.
		
		\bibitem{MR547117}
		W.~C. Waterhouse.
		\newblock {\em Introduction to affine group schemes}, volume~66 of {\em
			Graduate Texts in Mathematics}.
		\newblock Springer-Verlag, New York-Berlin, 1979.
		
	\end{thebibliography}
\end{document}